\documentclass[12pt]{amsart}
\usepackage{amsmath,amssymb}
\numberwithin{equation}{section}
\theoremstyle{plain}
  \newtheorem{theorem}{Theorem}[section]
  
  \newtheorem{proposition}[theorem]{Proposition}
  \newtheorem{lemma}[theorem]{Lemma}
  \newtheorem {corollary}[theorem]{Corollary}
 
\theoremstyle{definition}
  \newtheorem{definition}[theorem]{Definition}
  \newtheorem{example}[theorem]{Example}
\theoremstyle{remark}
 \newtheorem{remark}[theorem]{Remark}

\newcommand{\brac}[2]{{\langle #1\mid#2\rangle}}
\newcommand{\bri}[1]{{(#1)}}
\newcommand{\Aff}[1]{\widehat{#1}_{(\m,-\l)}}
\newcommand{\ab}[2]{{(a_{#1}^{(#2)},b_{#1}^{(#2)})}}
\def\bijec{\varrho}
\def\hom{\theta}
\def\homlm{\hom_{\lm,\mu}}

\def\dset{{\mathcal D}}
\def\con{{C}}
\def\cont{{\con_{\rbtab}}}
\def\wt{\zeta}
\def\tab{\mathop{\mathrm{Tab}}\nolimits}

\def\Tab{{\tab(\wh\lsm)}}
\def\TabR{{\tab^{\mathrm R}(\wh\lsm)}}
\def\TabRC{{\tab^{\mathrm{RC}}(\wh\lsm)}}
\def\V{\daff V}
\def\bx{\framebox(10,10)}
\def\proof{\noindent{\it Proof.}\ }
\def\aff{\dot}
\def\daff{\ddot}
\def\ch{^\vee}
\def\bra{{\langle}}
\def\ket{{\rangle}}
\def\qed{\hfill$\square$}
\def\udl{\underline}

\def\wh{\widehat}
\def\+{\mathop{\oplus}}
\def\*{\mathop{\otimes}} 
\def\pt{{u}}
\def\p{{p}}

\def\Q{{\mathbb Q}}
\def\Z{{\mathbb Z}}
\def\k{{\mathbb F}}
\def\F{{\mathbb Q}}

\def\W{{W}}
\def\affW{{\aff\W}}
\def\Wlm{{\aff W^{\lam-\mu}}}
\def\Zlm{{\aff Z^{\alsm}}}

\def\gl{{\mathfrak{gl}}}
\def\h{{\mathfrak h}}
\def\affh{{\aff\h}}
\def\tilh{{\tilde\h}}
\def\ii{{r}}
\def\al{{\alpha}}
\def\alch{\alpha^\vee}
\def\e{{\epsilon}}
\def\ech{\epsilon^\vee}
\def\lam{{\lambda}}
\def\lm{{\lambda}}
\def\lmm{{(\lambda,\mu)}}
\def\lsm{{\lambda/\mu}}
\def\alsm{{\widehat{\lam/\mu}}}

\def\ZZ{{\mathcal Z}}

\def\l{{\ell}}
\def\m{{m}}

\def\nxt{{{\tilde{\p}}}}
\def\t{{\tau}}

\def\Ical{{\widehat{\mathcal J}}^{n}}
\def\Icala{{\widehat{\mathcal J}}^{*n}}
\def\Icalml{{\Ical_{\m,\l}}}
\def\Icalaml{{\Icala_{\m,\l}}}
 
\def\affdom{\widehat{\mathcal P}^+}
\def\rbtab{T_{0}}

\def\X{{{\mathfrak X}}}

\def\cent{{\xi}}

\def\diagram{{{\Lambda}}}

\def\peri{{\gamma}}
\def\cat{{\mathcal O}^{ss}}
\def\Cat{{{{\mathcal O}_\kappa^{ss}(\daff H_n(q))}}}
\def\Sym{{\mathfrak S}}

\def\End{{\hbox{\rm{End}}}\,}
\def\Hom{{\hbox{\rm{Hom}}}\,}
\def\Irr{{\hbox{\rm{Irr}}}}
\def\dim{{\hbox{\rm{dim}}}\,}

\def\mod{{\rm mod\,}}
\def\gen{{{\rm{gen}}}}

\newcommand{\omitted}[1]{}
\begin{document}
\title[Double Affine Hecke Algebras]{Tableaux
on periodic skew diagrams and irreducible
representations of the double affine Hecke algebra of type $A$}
\author{
Takeshi Suzuki}
\address{Research Institute for Mathematical Sciences\\
Kyoto University, Kyoto, 606-8502, Japan}
\email{takeshi@kurims.kyoto-u.ac.jp}
\author{Monica Vazirani${}^*$}
\address{
 Department of Mathematics \\
UC Davis\\
One Shields Ave \\
Davis, CA 95616-8633}
\thanks{${}^*$Supported in part by an NSF
grant DMS-0301320.}
\email{vazirani@math.ucdavis.edu}

\begin{abstract}
We introduce and study an affine analogue of
 skew Young diagrams and tableaux on them.
 
It turns out that the double affine Hecke algebra 
of type $A$ acts on the space spanned by 
standard tableaux on each diagram. 
It is shown that the modules obtained this way are
 irreducible, and they exhaust all irreducible modules
 of a certain class over the double affine Hecke algebra.
In particular, the classification of 
 irreducible modules of this class, 
announced by Cherednik, is recovered.
\end{abstract}
\maketitle
\section{Introduction}
As is well-known, Young diagrams consisting of $n$ boxes
parameterize
isomorphism classes of
finite dimensional irreducible representations of the symmetric group 
$\Sym_n$ of degree $n$,
and moreover the structure of each irreducible representation
is described in terms of tableaux on the corresponding Young diagram;
namely, a basis of the representation is labeled by
 standard tableaux, with which
the action of $\Sym_n$ generators is explicitly described.
This combinatorial description due to A. Young
has played an essential role
in the study of the representation theory of 
the symmetric group (or the affine Hecke algebra),
and its generalization to 
the (degenerate) affine Hecke algebra $H_n(q)$ of $GL_n$
 has been given in \cite{Ch;special_bases,Ram1,Ram2},
where skew Young diagrams appear on 
combinatorial side. 

The purpose of this paper is to introduce
an ``affine analogue'' of skew Young diagrams and tableaux, which
 give a parameterization and a combinatorial description of
a family of irreducible representations of the double affine Hecke algebra 
$\daff H_n(q)$ 
of $GL_n$ over a field $\k$,
where $q\in\k$ is a parameter of the algebra.

The double affine Hecke algebra was introduced by 
I. Cherednik~\cite{Ch;unification,Ch;double} 
and has since been used by him and
 by several authors
to obtain important results
about diagonal coinvariants, Macdonald polynomials, and certain
Macdonald identities.

In this paper, we focus on the case where
$q$ is not a root of $1$,
and we consider representations of 
$\daff H_n(q)$ that are $\X$-semisimple;
namely, we consider representations which have 
basis of simultaneous eigenvectors with respect to
all elements in the commutative subalgebra
$\k[\X]=\k[x_1^{\pm1},x_2^{\pm1},\dots,x_n^{\pm1},\cent^{\pm1}]$
 of $\daff H_n(q)$.  
(In \cite{Ram1,Ram2}, such representations 
for affine Hecke algebras 
are referred to as ``calibrated'' )

On combinatorial side, we introduce
{\it periodic skew diagrams} as 
 skew Young diagrams consisting of
infinitely many boxes satisfying certain periodicity conditions. 
We define a  tableau on a periodic skew diagram
as a bijection from 
the diagram to $\Z$ 
which satisfies  the condition 
reflecting the periodicity of the diagram. 

Periodic skew diagrams 
are natural generalization of skew Young diagrams and 
have appeared 
in \cite{Ch;fourier} (or implicitly in \cite{AST}), but
the notion of tableaux on them seems new. 

To connect the combinatorics with the representation theory of 
the double affine Hecke algebra $\daff H_n(q)$,
we construct, for each periodic skew diagram,
an $\daff H_n(q)$-module that has a basis of
$\k[\X]$-weight vectors labeled by standard tableaux 
on the diagram
by giving the explicit action of the $\daff H_n(q)$ generators.

Such modules are $\X$-semisimple by definition.
We show that they are irreducible, and that 
our construction gives
a one-to-one 
 correspondence between
the set of periodic skew diagrams
and
the set of isomorphism classes of irreducible representations of 
the double affine Hecke algebra
 that are $\X$-semisimple.

The classification results here recover those of Cherednik's in
\cite{Ch;fourier} (see also \cite{Ch;introduction}), 
but in this paper we provide a detailed proof
based on purely combinatorial arguments  concerning 
standard tableaux on periodic skew diagrams.

Note that the corresponding
results for the degenerate double affine Hecke algebra  
of $GL_n$ easily follow from a parallel argument.


An outline of the paper is as follows.
Section\;\ref{section-affine.root} 
is a review
of the affine root system and 
the extended affine Weyl group of $\widehat{\gl}_n$.

The contents of Section\;\ref{section-periodic.skew} 
are purely combinatorial.
We introduce periodic skew diagrams and tableaux on them
in Section\;\ref{subsection-def.periodic.skew} and
Section\;\ref{subsection-def.tableaux} respectively.
These combinatorial objects  
are considered worth studying in themselves,
and here we investigate 
their relation with the affine Weyl group and 
{\it content\/}  
 functions.
The set of tableaux on a periodic skew 
diagram admits an action of the extended 
affine Weyl group $\aff W$, 
and it turns out that this action is simply transitive
and gives a bijective correspondence between 
the tableaux and the elements of $\aff W$.
In Section\;\ref{subsection;standard_tableaux},
we explicitly describe 
the subset $\aff W$
corresponding to the set of the standard tableaux,
which is the most interesting class 
from the view point of the representation theory.

We study 
content functions,
in particular 
those associated with standard tableaux,
in Section\;\ref{subsection-content.standard}.
The results obtained here lay the foundation to show 
our construction exhausts all $\X$-semisimple irreducible 
modules.

In Section\;\ref{section-rep.DAHA}, we 
introduce the double affine Hecke algebra
 and  
apply  the combinatorics 
studied in 
Section\;\ref{section-periodic.skew} 
to 
its  representation theory.

We remind
 the reader
in Section\;\ref{subsection-def.DAHA}
of the definition of the algebra
$\daff H_n(q)$ and review intertwining operators, 
which  were 
 also introduced by 
Cherednik and 
are elementary tools in the representation theory of 
$\daff H_n(q)$.

We derive some of rigid properties of 
$\X$-semisimple modules in Section\;\ref{section-Xssl}.
Then we give a combinatorial and explicit  construction
of  the representations of $\daff H_n(q)$
in Section\;\ref{ss;construction} using tableaux on 
periodic skew diagrams.
The related combinatorics is similar to and inspired by
that in \cite{Ram1} for the affine Hecke algebra.

Note that the statements of Section\;\ref{section-Xssl}
 also hold in the case that $q$ is a root of unity or for general $\kappa$,
where $\xi \in \daff H_n(q)$ is acting as the scalar $q^\kappa$.
However, when $q$ is a root of unity, the combinatorial description of the
modules is incredibly complicated.
 When $q$ is not a root of unity and we consider general,
not just integral $\kappa$, a combinatorial description of the
$\X$-semisimple modules can be reformulated from that described in the
paper.
The mathematics is essentially the same, but the way of drawing
and thinking about periodic skew diagrams becomes quite messy,
in contrast to the pretty pictures of 
Section\;\ref{section-periodic.skew}.

It is proved in Section\;\ref{subsection-classification}
 that we have constructed all the $\X$-semisimple 
irreducible representations
and that they are distinct up to diagonal shift 
of periodic skew diagrams.
This gives the classification of 
the $\X$-semisimple irreducible representations of $\daff H_n(q)$.

\medskip\noindent
{\bf Acknowledgment}

The second author would like to thank RIMS for their kind invitation
and hospitality during her visit.
\section{The affine root system and Weyl group}
\label{section-affine.root}
Let $\Q$ denote the field of rational numbers,
and let $\Z$ denote the ring of integral numbers.
We use the notation
\begin{equation*}
\Z_{\geq k}=\{n\in\Z\mid n\geq k\}
\end{equation*}
for $k\in\Z$, and
$$[i,j]=\{i,i+1,\dots,j\}$$
for $i,j\in\Z$ with $i\leq j$.

\subsection{The affine root system }
%

Let $n\in\Z_{\geq 2}.$
Let 
$\tilh$
be an $(n+2)$-dimensional vector space over $\F$
with the basis 
$\{\ech_1,\ech_2,\dots,\ech_n,c,d\}$:
$$\tilh=\left(\oplus_{i=1}^n \F\ech_i\right) \+\F c\+\F d.$$

Introduce the non-degenerate  symmetric bilinear form $(\ |\  )$ 
on $\tilh$ by
\begin{align*} 
(\ech_i|\ech_j)&=\delta_{ij},\quad 
(\e\ch_i| c)=(\e\ch_i|d)=0, \\
(c|d)&=1,\quad (c|c)=(d|d)=0.
\end{align*}
Put
$\h=\oplus_{i=1}^n \F\ech_i$
and $\affh=\h\+\F c$.

Let $\tilh^*=(\oplus_{i=1}^n\F\e_i)\oplus \F c^* \oplus \F \delta$ 
be the dual space of $\tilh$,
where $\e_i$,  $c^* $ and $\delta$
are   the dual vectors of $\ech_i$,
$c$ and $d$ respectively.

We  identify
the dual space $\affh^{*}$ 
of $\affh$ 
as a subspace 
of $\tilh^{*}$ via the identification
${\affh}^{*}={\tilh}^{*}/\F \delta\cong \h^{*}\oplus
\F c^{*}$.

The natural pairing is denoted by
$\bra\, |\,\ket : \tilh^*\times \tilh\to\F$. 
There exists an isomorphism
$\tilh^{*}\to\tilh$
such that 
$ \e_i\mapsto\ech_i$,
$\delta\mapsto c$ and $c^{*} \mapsto d$.
We denote by $\zeta\ch\in\tilh$ the image of
$\zeta\in\tilh^{*}$ under this isomorphism.
Introduce the bilinear form
$(\,\mid\, )$ on $\tilh^*$ through this 
isomorphism. Note that
$$(\zeta\mid \eta)=\brac{\zeta}{\eta\ch}=(\zeta\ch \mid \eta\ch),
\quad (\zeta,\eta\in\tilh^*).$$

Put $\al_{ij}=\e_i-\e_j\ (1\leq i\neq j\leq n)$ and 
$\al_i=\al_{ii+1}\ (1\leq i\leq n-1) $. 
Then 
\begin{align*}
R& =\left\{ \al_{ij}\mid i,j\in[1,n],\ i\ne j \right\},\
R^+=\left\{ \al_{ij}\mid i,j\in[1,n],\ i<j\right\},\\ 
\Pi&=\left\{\al_1,\al_2,\dots,\al_{n-1}\right\}
\end{align*}
give the system of roots,
positive roots and simple roots of type $A_{n-1}$
respectively.

Put $\al_0=-\al_{1n}+\delta$, and 
define the set $\aff R$ of (real) roots,  $\aff R^+$
of positive roots and  $\aff \Pi$ of simple roots 
of type $A_{n-1}^{(1)}$ by
$$
\begin{array}{l}
\aff R=
\left\{\al+k\delta\,|\,
\al\in R,\,k\in\Z\right\},\\
\aff R^+=\left\{\al+k\delta~|~
\al\in R^+,\,k\in \Z_{\geq0}\right\}
\sqcup \left\{-\al+k\delta~|~
\al\in R^+,\,k\in\Z_{\geq 1}\right\},\\
\aff\Pi=
\left\{\al_0,\al_1,\dots,\al_{n-1}
\right\}.
\end{array}
$$
\subsection{Affine Weyl group}
%
\begin{definition}
For $n\in\Z_{\geq 2}$,
the {\it extended affine Weyl group} $\affW_n$ of $\gl_n$
is the group defined by the following
generators and relations$:$
\begin{alignat*}{2}
&generators\ :&\ &  s_0,s_1,\dots,s_{n-1},
\pi^{\pm 1}.\\
&relations\ for\ n\geq 3
\  :&\ &  s_i^2=1\ (i\in[0,n-1]),\\
& & & s_i s_{j}\ s_i=
s_{j} s_i s_{j}\ (i-j\equiv\pm 1 \bmod n),\\
& & & s_i s_{j}=s_{j} s_i\ 
(i-j\not\equiv\pm 1 \bmod n),\\
& & & \pi s_i=s_{i+1} \pi,\ (i\in[0,n-2]),\quad
\pi s_{n-1}=s_0\pi,\\
& & & \pi\pi^{-1}=\pi^{-1}\pi=1.\\
&relations\ for\ n=2
\  :&\ &  s_0^2=s_1^2=1,\\
& & &\pi s_0=s_{1} \pi,\ \pi s_{1}=s_0\pi,\\
& & & \pi\pi^{-1}=\pi^{-1}\pi=1.
\end{alignat*}
\end{definition}
The subgroup $W_n$ of $\affW_n$
generated by
the elements 
$s_1,s_2,\dots,s_{n-1}$
is called the
{\it Weyl group} of $\gl_n$. 
The group $W_n$ is isomorphic to the symmetric group
of degree $n$.

In the following, we fix $n\in\Z_{\geq 2}$ and 
denote $\affW=\affW_n$ and $W=W_n$.

Put
$$P=\oplus_{i=1}^n \Z\e_i.
$$ 
Put $\t_{\e_1}=\pi s_{n-1}\cdots s_2s_1$ and
$\t_{\e_i}=\pi^{i-1}\t_{\e_1}\pi^{-i+1}\ (i\in[2,n])$.
Then there exists a group embedding $P\to \affW$
such that $\e_i\mapsto \t_{\e_i}$.
By $\t_\eta$ we denote the element in $\affW$
corresponding to $\eta\in P$.
It is well-known that the group $\affW$
is isomorphic to
the semidirect product $P\rtimes W$
with the relation $w \t_{ \eta} {w} ^{-1}=
\t_{{w}(\eta)}$.

The group $\affW$ acts on $\tilh$ by
\begin{alignat*}{1}
s_i (h)&=
 h-\bra\al_i| h\ket \alch_i\quad
 \hbox{for }i\in[1,n-1],\ h\in\tilh,\\
\t_{\e_i} ( h)&= h+
\bra\delta| h\ket \ech_i-
\left( \bra \e_i | h\ket+
\frac 12 
\bra\delta| h\ket\right)c\quad
 \hbox{for }i\in [1,n],\ h\in\tilh.
\end{alignat*}
The dual action on $\tilh^*$ is given by 
\begin{alignat*}{1}
s_{i}(\zeta)&=
\zeta-(\al_i|\zeta) \al_i\ \ 
\hbox{for }i\in[1,n-1],\ \zeta\in\tilh^*,\\
\t_{\e_i} (\zeta)&=
\zeta+(\delta|\zeta) \e_i-
\left( (\e_i|\zeta)+
\frac 12 (\delta| \zeta)
\right)\delta\quad
 \hbox{for }i\in [1,n],\ h\in\tilh^*.
\end{alignat*}
With respect to these actions, the inner products 
on $\tilh$ and $\tilh^*$ are $\affW$-invariant.
Note that the set $\aff R$ of roots is preserved by
the dual action of $\aff W$ on $\tilh^*$.
Note also that the action of $\aff W$ preserves
 the subspace  $\affh=\h\+\F c$, and
the dual action of $\affW$ on ${\affh}^*$ 
(called the affine action)
is described as follows:
\begin{equation}\label{eq;affine_action}
\begin{split}
s_i \left(\zeta\right)&=
\zeta-( \al_i| \zeta) \al_i\quad
\hbox{for }i\in[1,n-1],\ \zeta\in\affh^*,\\
\t_{\e_i} (\zeta)&=\zeta+(\delta|\zeta)  \e_i
\quad \hbox{for }i\in [1,n],\ h\in\affh^*.
\end{split}
\end{equation}
%
For $\al\in\aff R$, there exists 
$i\in[0,n-1]$ and $w\in \aff W$ such that
$w(\al_i)=\al$.
We set $s_\al=w s_i w^{-1}$.
Then
$s_\al$ is independent of the choice of $i$ and $w$, and we have
\begin{equation*}
  s_\al(h)=h-\brac{\al}{h}\alch
\end{equation*}
for $h\in\tilde{\h}$.
The element $s_\al$ is called the reflection corresponding to
$\al$. Note that $s_{\al_i}=s_i$.
%
For $w\in\affW $, set 
$$ R(w)=\aff R^+\cap w^{-1}\aff R^-,$$
where $\aff R^-=\aff R\setminus \aff R^+$.
The length $l(w)$ of $w\in\affW$ is defined as the number  
$\sharp R(w)$
of  elements in $ R(w)$.
For $w\in\affW $,
an expression $w=\pi^k s_{j_1}s_{j_2}\cdots s_{j_m}$ 
is called a reduced expression
if $m=l(w)$.
It can be seen that
\begin{equation}\label{eq;R(w)}
  R(w)=\{s_{j_m}\cdots s_{j_2}(\al_{j_1}),
s_{j_m}\cdots s_{j_3}(\al_{j_2}),\dots,\al_{j_m}\}
\end{equation}
if  $w=\pi^k s_{j_1} s_{j_2}\cdots s_{j_m}$ is a reduced 
expression.

Define the Bruhat order $\preceq$ in $\aff W$ by 
\begin{align*}
x\preceq w&\Leftrightarrow
x\hbox{ is equal to
a subexpression of a reduced expression of }w.
\end{align*}

We will review 
some fundamental facts 
in the theory of Coxeter groups, which are often used 
in this paper. See e.g. \cite{Hu} for proofs. 
\begin{lemma}
  \label{lem;coxeter}
$\rm{(i)}$
Let $w\in\affW$ and $i\in[0,n-1]$.
Then 
\begin{align*}
l(ws_i)&>l(w)\ \Leftrightarrow\ w(\al_i)\in \aff R^+,\\
l(s_iw)&>l(w)\ \Leftrightarrow\ w^{-1}(\al_i)\in \aff R^+.
\end{align*}

\smallskip\noindent
$\rm{(ii)}$ $($Strong exchange condition$)$
Let $\al\in\aff R^+$ and 
let $w\in\affW$ with a reduced expression 
$w=\pi^rs_{i_1}s_{i_2}\dots s_{i_k}$.
If $l(ws_\al)<l(w)$ then there exists $p\in[1,k]$
such that 
$ws_\al=\pi^rs_{i_1}s_{i_2}\dots\hat{s}_{i_p}\dots s_{i_k}$
$($omitting $s_{i_p})$.
Further $\al=s_{i_k}s_{i_{k+1}}\dots s_{i_{p+1}}(\al_{i_p})$.
\end{lemma}

Let $I$ be a subset of  $[0,n-1]$.
Put 
\begin{align*}
\aff\Pi_I&=\{\al_i\mid i\in I\}\subseteq\aff\Pi,\\
\affW_I&=\bra s_i\mid i\in I\ket\subseteq \affW,\\
\aff R_I^+&=\{\al\in\aff R^+\mid s_\al\in \affW_I\}.
\end{align*}
The subgroup $\affW_I$ is called the parabolic subgroup corresponding 
to $\aff\Pi_I$.
Define 
\begin{equation*}
 \affW^{I}=\left\{w\in \affW\mid R(w)\cap \aff R_I^+=\emptyset
\right\}.
\end{equation*}
The following fact is well-known.
\begin{proposition}
\label{pr;coset}
$\rm{(i)}$ $\aff W^I=\{w\in\aff W\mid
l(ws_{\al_i})>l(w)\text{ for all }\al_i\in\aff \Pi_I\}$.

\smallskip\noindent
$\rm{(ii)}$ 
For any $w\in\affW$, 
there exist a unique $x\in \affW^{I}$ 
and  a unique $y\in \affW_{I}$ such that $w=xy$.
Namely, 
the set $\affW^{I}$
gives a complete set of minimal length coset 
representatives for
$\affW/ \affW_{I}$.
\end{proposition}
\subsection{Notation}

For any integer $i$,
we introduce the following notation:
\begin{align}
 \e_i&=\e_{\udl{i}}-k\delta\in \tilh^*,\quad
\e\ch_i=\e\ch_{\udl i}-kc \in \tilh,
\end{align}
where  $i=\udl{i}+kn$ with
$\udl{i}\in[1,n]$ and $k\in \Z$.

Put 
$\al_{ij}=\e_i-\e_j$ 
and $\al\ch_{ij}=\e\ch_i-\e\ch_j$ for any $i, j\in\Z$.
Noting that $\e_0-\e_1=\delta+\e_n-\e_1=\al_0$, we 
reset $\al_{i}=\e_i-\e_{i+1}$ and $\al\ch_{i}=\e\ch_i-\e\ch_{i+1}$
for any $i\in \Z$.

The following is easy.
\begin{lemma}
\rm{(i)}  $\al_{i+n,j+n}=\al_{ij}$ for all $i,j\in\Z$. 

\smallskip\noindent
\rm{(ii)}  $\aff R=\{\al_{ij}\mid i,j\in\Z,\ 
i\not\equiv j \text{ mod }n \}$ as
a subset of $\tilh^*$.

\smallskip\noindent
\rm{(iii)} $\aff R^+=\{\al_{ij}\mid i,j\in\Z,\ i\not\equiv j \text{ mod }n
\text{ and } i< j\}$ as
a subset of $\tilh^*$.
\end{lemma}
%
%
Define the action of $\affW$ on the set $\Z$ 
of integers by
\begin{alignat*}{2}
s_i(j)&=j+1&\quad&\hbox{for }j\equiv i\ \mod n, \\
s_i(j)&=j-1&\quad&\hbox{for }j\equiv i+1\ \mod n,\\
s_i(j)&=j&\quad&\hbox{for }j\not\equiv i,i+1\ \mod n, \\
\pi(j)&=j+1&\quad&\hbox{for all }j.
\end{alignat*}
It is easy to see  that the action of $\t_{\e_i}$ $(i\in[1,n])$
is given by
\begin{alignat*}{2}
  \t_{\e_i}(j)&=j+n\quad& &\hbox{for }j\equiv i\ \mod n,\\
 \t_{\e_i}(j)&=j\quad& &\hbox{for }j\not\equiv i\ \mod n,
\end{alignat*}
and that the following formula holds for any $w\in\affW$:
$$  w(j+n)=w(j)+n\quad \hbox{for all }j.$$
\begin{lemma}\label{lem;w(e)}
Let $w\in\aff W$.

\smallskip\noindent
${\rm{(i)}}$
 $w(\e_j)=\e_{w(j)}$
and $w(\ech_j)=\ech_{w(j)}$ 
 for  any $j\in\Z$.

\smallskip\noindent
${\rm{(ii)}}$
 $w(\al_{ij})=\al_{w(i)w(j)}$
and $w(\alch_{ij})=\alch_{w(i)w(j)}$ 
 for  any $i,j\in\Z$.
\end{lemma}
\noindent {\it Proof.}
(i) 
It is enough to check the statement when $w=s_i~(i\in[1,n-1])$ 
and when $w=\t_{\e_i}~(i\in[1,n])$.
Let $j=\udl{j}+kn$ with $\udl{j}\in[1,n]$ and $k\in\Z$.
For $i\in[1,n-1]$, we have
$s_i(\e_j)=\e_j-(\al_i\mid\e_j)\al_i
=\e_j-(\al_i\mid\e_{\udl{j}})\al_i$.
This leads to $s_i(\e_j)=\e_{s_i(j)}$.
For $i\in[1,n]$, we have 
$\t_{\e_i}(\e_j)=\e_j+(\e_i\mid\e_j)\delta=
\e_j+(\e_i\mid\e_{\udl{j}})\delta$.
This leads to  $\t_{\e_i}(\e_j)=\e_{\t_{\e_i}(j)}$.

\noindent
(ii) The proof follows directly from (i).
\qed
\section{Periodic skew diagrams and tableaux on them}
\label{section-periodic.skew}
Throughout this paper, we let $\k$ denote
a field whose characteristic is not equal to $2$.
\subsection{Periodic skew diagrams}
\label{subsection-def.periodic.skew}
For $\m\in\Z_{\geq 1}$ and $\l\in\Z_{\geq 0}$, put
\begin{align}\label{eq;affdom}
\affdom_{\m,\l}
&=\{ \mu\in \Z^\m\mid
\mu_1\geq \mu_2\geq\cdots\geq\mu_\m\ \hbox{and }
\l\geq \mu_1-\mu_\m\},
\end{align}
where $\mu_i$ denotes the $i$-th component of $\mu$, i.e.,
$\mu=(\mu_1,\mu_2,\dots,\mu_\m)$.
Fix $n\in\Z_{\geq 2}$ and
introduce the following subsets of $\Z^\m\times \Z^\m$:
\begin{align*}
\Icalml&=
\left\{
\lmm\in \affdom_{\m,\l}\times \affdom_{\m,\l}\left|
\lam_i\geq \mu_i\ (i\in[1,\m]),\
\sum\limits_{i=1}^\m(\lam_i-\mu_i)=n 
\right.\right\},\\
\Icalaml&=
\left\{
\lmm\in \affdom_{\m,\l}\times\affdom_{\m,\l} \left|
\lam_i>\mu_i\ (i\in[1,\m]),\
\sum\limits_{i=1}^\m(\lam_i-\mu_i)=n 
\right.\right\}.
\end{align*}
For $\lmm\in \Icalml$, 
define
the subsets $\lsm$ and   $\Aff{\lsm}$ 
of $\Z^2$ by
\begin{align*}
  \lsm&=\left\{(a,b)\in\Z^2 \mid
a\in [1,\m],\ b\in [\mu_a+1,\lam_a] \right\},\\
 \Aff{\lsm}&=\left\{(a+k\m,b-k\l)\in\Z^2 \mid
(a,b)\in\lsm,\ k\in\Z \right\}.
\end{align*}
Let $\lsm[k] = \lsm + k(\m,-\l)$. 
Obviously we have
$$\Aff{\lsm}=\bigsqcup_{k\in\Z} \lsm [k] =
\bigsqcup_{k\in\Z}\left(\lsm+k(\m,-\l)\right).$$ 

The set  $\lsm$ is
so called 
the skew diagram (or skew Young diagram) 
associated with $\lmm$.

We call the set  $\Aff\lsm$ 
the periodic skew diagram
associated with $\lmm$.

We will denote $\Aff\lsm$ just by $\alsm$ 
when 
 $\m$ and $\l$ are fixed.
\begin{example}\label{ex1}
(i) Let $n=7,~\m=2$ and $\l=3$.
Put $\lm=(5,3),~\mu=(1,0)$.
Then $\lmm\in\Icalaml$ and we have
\begin{equation*}
  \lsm=\{(1,2),(1,3),(1,4),(1,5),(2,1),(2,2),(2,3)\}.
\end{equation*}
The set $\lsm$ is expressed by the following picture
(usually, the coordinate in the boxes are omitted):

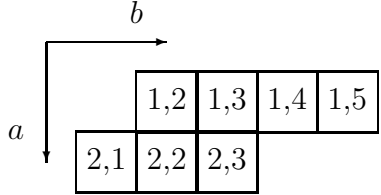
\begin{center}
\setlength{\unitlength}{.08cm}
\begin{picture}(60,40)(0,-15)
\put(10,0){\bx{$1,\!2$}}
\put(20,0){\bx{$1,\!3$}}
\put(30,0){\bx{$1,\!4$}}
\put(40,0){\bx{$1,\!5$}}
\put(0,-10){\bx{$2,\!1$}}
\put(10,-10){\bx{$2,\!2$}}
\put(20,-10){\bx{$2,\!3$}}
\put(-5,15){\vector(1,0){20}}
\put(-5,15){\vector(0,-1){20}}
\put(5,15){\makebox(10,10){$b$}}
\put(-15,-5){\makebox(10,10){$a$}}
\end{picture}
\end{center}
The periodic skew diagram
$$\alsm_{(2,-3)}= \bigsqcup_{k\in\Z} \lsm [k]
=\bigsqcup_{k\in\Z}\left( \lsm+k(2,-3)\right)$$
is expressed by the following picture:
%
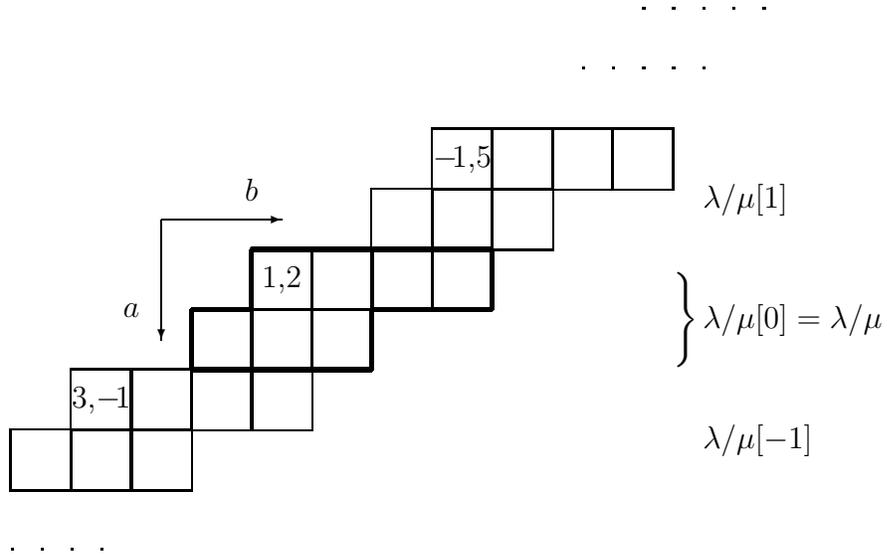
\begin{center}
\setlength{\unitlength}{.08cm}
\begin{picture}(150,105)(-40,-45)
\multiput(75,50)(5,0){5}{\rule{.8pt}{.8pt}}
\multiput(65,40)(5,0){5}{\rule{.8pt}{.8pt}}
\put(40,20){\bx{$-\!1,\! 5$}}
\multiput(40,20)(10,0){4}{\bx}
\multiput(30,10)(10,0){3}{\bx}
\put(85,17){$\lsm [1]$}
\put(10,0){\bx{$1,\!2$}}
\multiput(10,0)(10,0){4}{\bx}
\multiput(0,-10)(10,0){3}{\bx}
\put(85,-3){$\lsm [0] = \lsm$} 
\put(-20,-20){\bx{$3,\! -\!1$}}
\multiput(-20,-20)(10,0){4}{\bx}
\multiput(-30,-30)(10,0){3}{\bx}
\put(85,-23){$\lsm [-1]$}
\multiput(-35,-40)(5,0){5}{\rule{.8pt}{.8pt}}
\put(-5,15){\vector(1,0){20}}
\put(-5,15){\vector(0,-1){20}}
\put(5,15){\makebox(10,10){$b$}}
\put(-15,-5){\makebox(10,10){$a$}}
\put(80,-3){\mbox{$\Biggl. \Biggr\}$}}
\linethickness{1.5pt}
\put(10,10){\line(1,0){40}}
\put(10,10){\line(0,-1){10}}
\put(0,0){\line(1,0){10}}
\put(0,0){\line(0,-1){10}}
\put(0,-10){\line(1,0){30}}
\put(30,0){\line(1,0){20}}
\put(30,0){\line(0,-1){10}}
\put(50,0){\line(0,1){10}}
\end{picture}
\end{center}
%
\end{example}
Generally, periodic skew diagrams are defined as follows 
(see \cite{Ch;fourier}):
\begin{definition}
\label{def;periodicdiagram}
For $\peri\in\Z^2$,
a subset $\diagram\subset \Z^2$ 
is called
 a $\peri$-{\it periodic skew diagram} 
(or a {\it periodic skew diagram} of {\it period} $\peri$)
if it satisfies the following conditions. 

\smallskip\noindent
${\rm (D1)}$ 
The set $\diagram$ is invariant
 under the parallel translation by $\peri:$
$$\diagram+\peri=\diagram,$$
and hence the group $\Z \peri$ acts on  $\diagram$.

\smallskip
\noindent
${\rm (D2)}$ 
A fundamental domain of the action of $\Z\peri$
on $\diagram$ consists of finitely many elements.
This number is 
called the {\it degree} of $\diagram$.

\smallskip
\noindent
${\rm (D3)}$
If $(a,b)\in \diagram$ and $(a+i,b+j)\in\diagram$ for 
$i,j\in\Z_{\geq 0}$,
then  
the rectangle
$\{(a+i',b+j')\mid i'\in [0,i],\ j'\in[0,j]\}$
is included in $\diagram$.
%
\end{definition}
Let $ \dset^n_\peri$ denote
the set of all $\peri$-periodic skew diagram of degree $n$, and put
\begin{align*}
\dset^{*n}_\peri
&=\{\diagram\in\dset^n_\peri\mid
\forall a\in\Z,
\exists b\in\Z\text{ such that }
(a,b)\in\diagram\}.
\end{align*}
Namely,  $\dset^{*n}_\peri$ is
 the subset of $\dset^n_\peri$
consisting of
all diagrams without empty rows.

Note that an element in $\dset_{(0,0)}^n$ is 
regarded as a (classical) skew Young diagram of degree $n$.
 \begin{lemma}
   Let $\peri\in\Z^2$.

\smallskip\noindent
${\rm(i)}$
If $\dset^{n}_\peri\neq \emptyset$,
then $\peri\in \Z_{\leq 0}\times \Z_{\geq 0}$
or $\peri\in \Z_{\geq 0}\times \Z_{\leq 0}$.

\smallskip\noindent
${\rm(ii)}$
If $\dset^{*n}_\peri\neq \emptyset$,
then $\peri\in \Z_{\geq 1}\times \Z_{\leq 0}$
or $\peri\in \Z_{\leq -1}\times \Z_{\geq 0}$.
 \end{lemma}
\proof
(i) Since $\dset^{n}_\peri=\dset^{n}_{-\peri}$,
it is enough to prove that
$\dset^{n}_\peri=\emptyset$ for $\peri\in\Z_{\geq 1}\times\Z_{\geq 1}$.

Suppose $\dset^{n}_{(\m,\l)}\neq\emptyset$
for some $\m\in\Z_{\geq 1}$ and $\l\in\Z_{\geq 1}$, and
take $\diagram\in\dset^{n}_{(\m,\l)}$.
Then $\bar{\diagram}=\{(a,b)\in\diagram\mid a\in[1,\m]\}$
is a  fundamental domain of the action of $\Z (\m,\l)$ on $\diagram$.

Let $(a,b)\in \bar\diagram$.
Then the condition (D1) implies
$\{(a+k\m,b+k\l)\}_{k\in\Z}\subseteq \diagram$,
and the
condition (D3) implies 
$\{(a,b+k\l)\}_{k\in\Z}\subseteq\diagram$
and hence $\{(a,b+k\l)\}_{k\in\Z}\subseteq\bar\diagram$.
This implies that the fundamental domain $\bar{\diagram}$
contains infinitely many elements.
This contradicts the condition (D2), and 
hence we have
 $\dset^{n}_{(\m,\l)}=\emptyset$ for  
$\m\in\Z_{\geq 1}$ and $\l\in\Z_{\geq 1}$. 

\noindent
(ii) By (i), it is enough to show that 
 $\dset^{*n}_{(0,\l)}=\emptyset$ for all $\l\in\Z$,
and this is easy.\qed

\medskip
Let $\m\in\Z_{\geq 1}$ and $\l\in\Z_{\geq 0}$.
For $\lmm\in\Icalml$, it is easy to see that
the set $\Aff\lsm$ satisfies the conditions (D1)(D2)(D3) 
in Definition~\ref{def;periodicdiagram}, and hence we have
$\Aff\lsm\in\dset^{n}_{(\m,-\l)}$.
%
%
%
\begin{proposition}\label{pr;affine_diagram}
Let $n\in\Z_{\geq 2},\ 
\m\in\Z_{\geq 1}$ and $\l\in\Z_{\geq0}$.
The correspondence 
$\Icalml{\rightarrow}\dset^{n}_{(\m,-\l)}$
given by $\lmm\mapsto\alsm$
is a surjection. %
Moreover, its restriction to $\Icalaml$ gives
 a bijection
$$\Icalaml\overset{\sim}{\rightarrow}\dset^{*n}_{(\m,-\l)}.$$
\end{proposition}
\noindent
{\it Proof.}
Take any  $\diagram\in\dset^{n}_{(\m,-\l)}$.

Fix $i_0\in\Z_{\leq 0}$ such that the $i_0$-th row
of $\diagram$ is not empty.
For  $i\geq i_0$, define $\lm_i$ and $\mu_i$ 
recursively by the following relations:
\begin{align*}
  \lm_i&=
\begin{cases}
\text{max}\{b\in\Z\mid (i,b)\in\diagram\}
&\text{ if } i\hbox{-th row is not empty},\\
\lm_{i-1}
&\text{ if } i\hbox{-th row is empty},
\end{cases}\\
 \mu_i&=
\begin{cases}
\text{min}\{b\in\Z\mid (i,b)\in\diagram\}-1
&\text{ if } i\hbox{-th row is not empty},\\
\lm_{i-1}
&\text{ if } i\hbox{-th row is empty}.
\end{cases}
\end{align*}
Put $\lm=(\lm_1,\lm_2,\dots,\lm_\m)$ and 
$\mu=(\mu_1,\mu_2,\dots,\mu_\m)$.
Then it follows from the condition (D3) (with $i=0$)
that  $$\{(a,b)\in\diagram\mid a=i\}=[\mu_i+1,\lm_i]
\quad (i\in[1,\m])$$ and hence
$$
\lsm=\{(a,b)\in\diagram\mid
a\in[1,\m]\}.$$
It follows from the condition (D1) 
that
$\lsm$ is a fundamental domain of $\Z(\m,-\l)$ on $\diagram$
and
\begin{equation}
\label{eq;diagram_sum}
\diagram=\bigsqcup_{k\in\Z}
(\lsm+k(\m,-\l))=\alsm_{(\m,-\l)}.
\end{equation}
In particular, we have $\sharp\lsm=n$.

Note that 
$\lm_{0}=\lm_\m+\l$ and $\mu_{0}=\mu_\m+\l$
by the condition (D1).

Now, the condition (D3) implies that $\lm_i\geq \lm_{i+1}$
and $\mu_i\geq\mu_{i+1}$ for all $i\geq i_0$, in particular,
for all $i\in[1,m-1]$.
This yields
 $\lm\in\affdom_{\m,\l}$ and $\mu\in\affdom_{\m,\l}$.
Therefore the correspondence $\Icalml\to\dset^n_{(\m,-\l)}$
is surjective.

Now, it is clear from the discussion above
that the correspondence $\lmm\mapsto\alsm$
gives a bijection
$\Icalaml\to\dset^{*n}_{(\m,-\l)}$.
\qed
\subsection{Tableaux on periodic skew diagram}
\label{subsection-def.tableaux}
Fix $n\in\Z_{\geq 2}$.
Recall that a  bijection from
a skew Young diagram, say $\lsm$, of degree $n$ 
to the set $[1,n]$ is
called a tableau on $\lsm$.
\begin{definition}\label{def;tableau}
For $\peri\in\Z^2$ and $\diagram\in\dset_\peri^n$,
 a bijection
$T:\diagram\to\Z$
is said to be a $\peri$-{\it tableau} on $\diagram$
if $T$ satisfies 
\begin{equation}
  \label{eq;Tab}
T(\pt+\peri)=T(\pt)+n \quad\text{for all }\pt\in\diagram.
  \end{equation}
\end{definition}
Let $\tab_{\peri}(\diagram)$ denote 
the set of all $\peri$-tableaux on $\diagram$.

In this paper, we mostly treat periodic skew diagrams associated with 
$\lmm\in\Icalml$ for some $\m\in\Z_{\geq 1}$ and $\l\in\Z_{\geq 0}$.
For $\lmm\in\Icalml$, we always 
choose $(\m,-\l)$ as a period of $\alsm$.
We use 
 the abbreviated notation
$$\tab(\alsm)=\tab_{(\m,-\l)}(\alsm)$$ for $\lmm\in\Icalml$,
and we let a tableau on $\alsm$ mean  
a $(\m,-\l)$-tableau on $\alsm$.
\begin{remark}
A tableau on $\alsm$ is determined uniquely 
from the values on  a fundamental domain of $\alsm$ 
with respect to the action of  $\Z\peri$.
It also holds that any bijection
from a fundamental domain of $\Z\peri$ to the set $[1,n]$  
uniquely extends to a tableau on $\alsm$.
\end{remark}
There exists a unique tableau
$\rbtab^{\alsm}={\rbtab}$ on $\alsm$ such that
\begin{equation}\label{eq;rrtab}
\rbtab (i,\mu_i+j)=\sum_{k=1}^{i-1}(\lm_k-\mu_k)+j\quad
\text{for}\ i\in[1,\m],\ j\in[1,\lm_i-\mu_i].
\end{equation}
We call $\rbtab$ the {\it row reading tableau} on $\alsm$.
\begin{example}
Let $n=7,~\m=2$, $\l=3$
and $\lm=(5,3),~\mu=(1,0)$.
The 
tableau $\rbtab$ on $\alsm$ given above
is expressed as follows:
%
\begin{center}
\setlength{\unitlength}{.08cm}
\begin{picture}(150,90)(-40,-45)
\multiput(65,40)(5,0){5}{\rule{.8pt}{.8pt}}
\put(40,20){\bx{$-\!6$}}
\put(50,20){\bx{$-\!5$}}
\put(60,20){\bx{$-\!4$}}
\put(70,20){\bx{$-\!3$}}
\put(30,10){\bx{$-\!2$}}
\put(40,10){\bx{$-\!1$}}
\put(50,10){\bx{$0$}}
\put(10,0){\bx{$1$}}
\put(20,0){\bx{$2$}}
\put(30,0){\bx{$3$}}
\put(40,0){\bx{$4$}}
\put(0,-10){\bx{$5$}}
\put(10,-10){\bx{$6$}}
\put(20,-10){\bx{$7$}}
\put(85,-3){$\lsm$}
\put(-20,-20){\bx{$8$}}
\put(-10,-20){\bx{$9$}}
\put(0,-20){\bx{$10$}}
\put(10,-20){\bx{$11$}}
\put(-30,-30){\bx{$12$}}
\put(-20,-30){\bx{$13$}}
\put(-10,-30){\bx{$14$}}
\multiput(-35,-40)(5,0){5}{\rule{.8pt}{.8pt}}
\put(80,-3){\mbox{$\Biggl. \Biggr\}$}}
\linethickness{1.5pt}
\put(10,10){\line(1,0){40}}
\put(10,10){\line(0,-1){10}}
\put(0,0){\line(1,0){10}}
\put(0,0){\line(0,-1){10}}
\put(0,-10){\line(1,0){30}}
\put(30,0){\line(1,0){20}}
\put(30,0){\line(0,-1){10}}
\put(50,0){\line(0,1){10}}
\end{picture}
\end{center}
\end{example}
\begin{proposition}
  Let $\lmm\in\Icalml$.
The group $\affW$ acts on the set $\Tab$ by
  \begin{equation}\label{eq;actonTab}
    (wT)(\pt)=w(T(\pt))
  \end{equation}
for $w\in \affW$, $T\in\Tab$ and $\pt\in\wh\lsm$.
\end{proposition}
\noindent {\it Proof.}
It is obvious that $wT$ is a bijection.
It is enough to verify that $wT$ 
satisfies the condition 
\eqref{eq;Tab} in Definition~\ref{def;tableau}.
Putting $\peri=(\m,-\l)$, we have
$$
 (wT)(\pt+\peri)=w(T(\pt+\peri))
=w(T(\pt)+n)=wT(\pt)+n.
$$
Therefore $wT$ satisfies \eqref{eq;Tab}.
\qed

\medskip
For each  $T\in\Tab$, 
define the map
$$\psi_T:\affW\to\Tab$$
by $\psi_T(w)= wT$ $(w\in\affW)$.
\begin{proposition}\label{pr;W=Tab}
Let $\lmm\in\Icalml$.
For any $T\in\Tab$, the correspondence 
$\psi_T$ is a bijection.
\end{proposition}
\noindent {\it Proof.}
It is enough to show the statement for $T={\rbtab}$
given by \eqref{eq;rrtab}.
We prove the surjectivity first.
Take any $S\in\Tab$ and
put $\ii_i=S(\rbtab^{-1}(i))$ for $i\in[1,n]$. 
Suppose $\ii_i-\ii_j=kn$ for some  $i,j\in[1,n]$ and some $k\in\Z$.
Then $S(\rbtab^{-1}(j)+k(\m,-\l))=\ii_j+kn=\ii_i
=S(\rbtab^{-1}(i))$,
and hence
$\rbtab^{-1}(j)+k(\m,-\l)= \rbtab^{-1}(i)$.
This means $k=0$ and $i=j$.

Let $\ii_i=\udl{\ii}_i+k_in$
with $\udl{\ii}_i\in[1,n] $ and $k_i\in\Z$.
Then we have shown that
$\udl{\ii}_i\neq\udl{\ii}_j$ for $i,j\in[1,n]$ such that $i\neq j$.
This ensures that there exists $x\in\W$ such that 
 $x(i)=\udl{\ii}_i$ for all $i\in[1,n]$.
Putting $w:=x\cdot \t_{\e_1}^{k_1} \t_{\e_2}^{k_2}\cdots
\t_{\e_n}^{k_n}$, we have 
$w(i)=\ii_i$ for any $i\in[1,n]$ and hence $w{\rbtab}=S$
on the fundamental domain $\lsm$.
This implies $w\rbtab=S$.

It is easy to see that the choice of $w$ for each $S$ is unique,
and hence the injectivity follows.
\qed\par\medskip
The following formula follows directly from 
the definition 
\eqref{eq;actonTab} 
of the action of $\affW$.
\begin{lemma}\label{lem;Twinv}
  $T^{-1}(w^{-1}(i))=(wT)^{-1}(i)$ for any  $T\in\Tab$,  $w\in\affW$
and $i\in\Z$.
\end{lemma}
\subsection{Content and weight}
\label{subsection-content.weight}
Let $\con$ denote the map from $\Z^2$ to $\Z$
given by $\con(a,b)=b-a$ for $(a,b)\in\Z^2$.

For a tableau $T\in\Tab$, 
define the map $\con^{\alsm}_T:\Z\to\Z$ 
by
$$\con^\alsm_T(i)=\con(T^{-1}(i))\quad (i\in\Z),$$
and call $\con^\alsm_T$ the {\it content} 
of $T$.
We simply denote $\con^\alsm_T$ by $\con_T$ when
$\lmm$ is fixed.
\begin{lemma}\label{lem;content}
Let $T\in\Tab$. Then

\smallskip\noindent
${\rm{(i)}}$
$  \con_T(i+n)=\con_T(i)-(\l+\m)$ for all $i\in\Z$.

\smallskip\noindent
${\rm{(ii)}}$
$\con_{wT}(i)=\con_T(w^{-1}(i))$ for all $w\in\aff W$ and $i\in\Z$.
\end{lemma}
\noindent
{\it Proof.}
(i) Put $(a,b)=T^{-1}(i)\in \alsm$.
Then $T(a+\m,b-\l)=T(a,b)+n=i+n$.
We have
\begin{align*}
\con_T(i+n)&=\con(T^{-1}(i+n))=\con(a+\m,b-\l)\\
&=(b-a)-(\l+\m)=\con_T(i)-(\l+\m).
\end{align*}
(ii) The proof directly follows from Lemma~\ref{lem;Twinv}. 
\qed

\medskip
For $T\in\Tab$, we define
$\wt_T\in \affh^*$ by
$$\wt_T=\sum_{i=1}^n\con_T(i)\e_i+(\l+\m)c^*.$$
Then $\wt_T$ belongs to the lattice
$$\aff P\overset{def}{=}P\oplus\Z c^*=
(\mathop{\oplus}_{i=1}^n\Z\e_i)\oplus\Z c^*.$$
Note that the action \eqref{eq;affine_action}
of $\aff W$ on $\affh^*$ preserves $\aff P$.
Lemma~\ref{lem;content} immediately implies the following:
\begin{lemma}
Let
$T\in\Tab$. Then

\smallskip\noindent
$\rm{(i)}$
$\brac{\wt_T}{\ech_i}={\con_T(i)}\hbox{ for all }i\in\Z.$

\smallskip\noindent
$\rm{(ii)}$
$w(\wt_T)=\wt_{wT} \hbox{ for all }w\in\aff W.$
\end{lemma}
%
\subsection{The affine Weyl group and row increasing tableaux}
\label{subsection;row_increasing_tableaux}
Let $\lmm\in\Icalml$.
\begin{definition}
A tableau $T\in\Tab$ is said to be {\it row increasing}
(resp. {\it column increasing})
if 
\begin{equation*}
\begin{split}
  (a,b), (a,b+1)\in\alsm\ &\Rightarrow\ T(a,b)<T(a,b+1).\\
(\hbox{resp. }
  (a,b), (a+1,b)\in\alsm\ 
&\Rightarrow\ T(a,b)<T(a+1,b).
)
\end{split}
\end{equation*}
A tableau $T\in\Tab$ which 
is
row increasing and column increasing is called 
{\it a standard tableau}
(or {\it a row-column increasing tableau}).

Denote by $\TabR$ (resp.\;$\TabRC$)
the set of all row increasing (resp. standard)
tableaux on $\alsm$.
\end{definition}

For $\lmm\in\Icalml$,
put
$$I_{\lm,\mu}=[1,n-1]\setminus\{n_1,n_2,\dots,n_{\m-1}\},$$
where
$n_i=\sum_{j=1}^i(\lm_j-\mu_j)$
for $i\in[1,\m-1].$

We write  $\aff R_{\lm-\mu}^+=\aff R_{I_{\lm,\mu}}^+$,
$\aff W_{\lm-\mu}=\affW_{I_{\lm,\mu}}$ 
and
$\Wlm=\affW^{I_{\lm,\mu}}$.

Note that  $\aff R_{\lm-\mu}^+\subseteq R^+$ and 
$\affW_{\lm-\mu}=\W_{\lm_1-\mu_1}\times
\W_{\lm_2-\mu_2}\times\cdots\times\W_{\lm_\m-\mu_\m}
\subseteq \W$.

Recall that  the correspondence
 $\psi_T:\affW\to \Tab$ given by $w\mapsto w{T}$
is bijective (Proposition~\ref{pr;W=Tab})
for any $T\in\Tab$.
\begin{proposition}\label{pr;Wlm=TabR}
Let $\lmm\in\Icalml$. Then 
$$\psi_{T_0}^{-1}(\TabR)=\Wlm,$$
or equivalently, $\TabR=\Wlm \rbtab
{=}\{w\rbtab\mid w\in\Wlm\}.$
\end{proposition}
\noindent {\it Proof.}
First 
we will prove $\Wlm \rbtab \subseteq \TabR$.

Take
 $(a,b),(a,b+1)\in\alsm$ and put ${\rbtab}(a,b)=i$.
Then ${\rbtab}(a,b+1)=i+1$
and $\al_i\in \aff R_{\lm-\mu}^+$.
If $w\in\Wlm$, then we have
 $l(ws_i)>l(w)$.
This means $w(\al_i)=\e_{w(i)}-\e_{w(i+1)}\in \aff R^+$.
Hence $w(i)<w(i+1)$, or equivalently $w{\rbtab}(a,b)<w{\rbtab}(a,b+1)$.
Therefore $w{\rbtab}\in\TabR$
 for all $w\in\Wlm$. 

Next, we will prove $\Wlm T_0\supseteq\TabR$.

For $T\in\TabR$, take $w\in\affW$ such that $w{\rbtab}=T$.
We have to show that $w\in \affW^{\lm-\mu}$.
Let $\al_{ij}\in \aff R_{\lm-\mu}^+$.
Put $(a,b)={\rbtab}^{-1}(i).$
Then ${\rbtab}^{-1}(j)=(a,b+j-i).$
Since $w{\rbtab}$ is row increasing, we have
$w{\rbtab}(a,b)<w{\rbtab}(a,b+j-i)$
and hence $w(i)<w(j)$.
This means that $w(\al_{ij})\in\aff R^+.$
Therefore 
$\al_{ij}\not\in
R(w)=\aff R^+\cap w^{-1}\aff R^-$.
This proves $R(w)\cap \aff R_{\lm-\mu}^+=\emptyset$
and hence $w\in\Wlm$.
\qed
\subsection{The set of standard tableaux}
\label{subsection;standard_tableaux}

The next lemma follows easily:
\begin{lemma}\label{lem;diagonalinc}
Let $\lmm\in\Icalml$ and $T\in\TabRC$.
If $(a,b)\in\alsm$ and  $(a+1,b+1)\in\alsm$,
then $T(a+1,b+1)-T(a,b)>1$.
\end{lemma}
As a direct consequence of
Lemma\ref{lem;diagonalinc}, we obtain
the following result, which will be used in the next section:
%
\begin{proposition}\label{pr;CT=CS}
  Let $\lmm\in\Icalml$ and  $T,S\in\TabRC$.
If
$\con_T=\con_S$ then
$T=S.$
\end{proposition}
%
Our next purpose is to describe
the subset of $\affW$ which corresponds to
$\TabRC$ under the correspondence $\psi_T$
$(T\in\TabRC)$.
\begin{lemma}\label{lem;wBRC}
  Let $\lmm\in\Icalml$.
Let $w\in\affW$ and $i\in[0,n-1]$ such that
 $w{\rbtab}\in\TabRC$ and $l(w)>l(s_iw)$. Then $s_iw{\rbtab}\in\TabRC$.
\end{lemma}
\noindent {\it Proof.}
We have $w\in\Wlm$ by Proposition~\ref{pr;Wlm=TabR}.
Put $x=s_iw$.
Since $R(w)=R(x)\sqcup\{x^{-1}(\al_i)\}$,
we have  $x\in\Wlm$.
Hence $x{\rbtab}\in\TabR$.
  
Suppose $x{\rbtab}\not\in\TabRC$.
Then there exist $(a,b)$ and $(a+1,b)$ in $\alsm$ 
such that
\begin{equation*}
\begin{split}
  x{\rbtab}(a,b)&>x{\rbtab}(a+1,b),\\
  w{\rbtab}(a,b)&<w{\rbtab}(a+1,b).
\end{split}
\end{equation*}
This implies $x{\rbtab}(a,b)=i+1+kn$
and $x{\rbtab}(a+1,b)=i+kn$ for some $k\in\Z$.

On the other hand, we have $x^{-1}(\al_i)\in \aff R^+$ 
as $l(w)>l(x)$.
Therefore it follows that $x^{-1}(i)<x^{-1}(i+1)$
and hence  $x^{-1}(i+kn)<x^{-1}(i+1+kn)$.
Therefore we have
${\rbtab}(a+1,b)<{\rbtab}(a,b)$, and
this is a contradiction.
\qed

\medskip
For $T\in\TabRC$, put
\begin{align}
\label{eq;Zlm}
&\Zlm_{T}=
\left\{\left.
w\in\affW\right|
\brac{\wt_T}{\al\ch}
\not\in \{-1,1\}\hbox{ for all }
\al\in R(w)\right\}.
\end{align}
\begin{lemma}
  \label{lem;ZlmsubWlm}
Let $\lmm\in\Icalml$. Then 
$$\Zlm_{\rbtab}\subseteq\Wlm.$$
\end{lemma}
\proof
Take $w\in\affW$ such that $w\notin\Wlm$.
Then, it follows that there exists $j\in[0,n-1]$ such that
 $s_j\in \affW_{\lm-\mu}$ and $l(ws_j)<l(w)$.
Then Lemma~\ref{lem;coxeter}-(ii) implies that $\al_j  \in R(w)$.
By $s_j\in \affW_{\lm-\mu}$,
we have $\brac{\wt_{\rbtab}}{\alch_j} =-1$.
Hence we have $w\not\in\Zlm_{\rbtab}$, and
proved $\Zlm_{\rbtab}\subseteq \Wlm$.
\qed
%
\begin{theorem}
\label{th;Zlm=TabRC2}
Let $\lmm\in\Icalml$ and $T\in\TabRC$. Then 
$$\psi_T^{-1}(\TabRC)=\Zlm_T,$$
or equivalently, 
$\TabRC=\Zlm_T T.$  
\end{theorem}
\noindent {\it Proof.}    
({\it Step 1})
First we will prove the statement for the row reading tableau
$\rbtab$, namely, we will prove $\TabRC=\Zlm_{\rbtab} \rbtab.$

Let us see $\Zlm_{\rbtab}\rbtab\subseteq\TabRC$,
that is,  $w{\rbtab}\in\TabRC$ for all $w\in\Zlm_{\rbtab}$.
We proceed by induction on $l(w)$.

If $w\in\Zlm_{\rbtab}$ with $l(w)=0$ then $w=\pi^k$ for some $k\in\Z$
and  it is obvious that $w{\rbtab}$ is row-column increasing.
Suppose that $w{\rbtab}$ is row-column increasing
for all $w\in\Zlm_{\rbtab}$ with $l(w)<k$.

Take $w\in\Zlm_{\rbtab}$ with $l(w)=k$. 
Note that $w\in\Wlm$ by Lemma~\ref{lem;ZlmsubWlm}.
Take $x\in\aff W$ and $i\in[0,n-1]$ such that
$w=s_ix$ and $l(w)=l(x)+1$.

Note that 
$R(w)=R(x)\sqcup\{x^{-1}(\al_i)\}$,
and hence 
$x\in\Zlm_{\rbtab}$.
By the induction hypothesis we have $x{\rbtab}\in \TabRC$. 
Suppose that $w{\rbtab}\not\in \TabRC$. 
Then $w{\rbtab}$ is not column increasing 
because $w{\rbtab}$ is row increasing by Proposition~\ref{pr;Wlm=TabR}.
Therefore there exist $(a,b),(a+1,b)\in\alsm$
such that
\begin{equation*}
  \begin{split}
x{\rbtab}(a,b)&<x{\rbtab}(a+1,b),\\
w{\rbtab}(a,b)&>w{\rbtab}(a+1,b).
  \end{split}
\end{equation*}
This implies that 
$x{\rbtab}(a,b)=i+kn$ and $x{\rbtab}(a+1,b)=i+1+kn$
for some $k\in\Z$.
We have
\begin{equation*}
\begin{split}
\brac{\wt_{\rbtab}}{x^{-1}(\alch_i)}
=&\brac{\wt_{\rbtab}}{x^{-1}(\alch_{i+kn})}\\
=&\cont(x^{-1}(i+kn))-\cont(x^{-1}(i+1+kn))\\
=&b-a-(b-(a+1))=1.
\end{split}
\end{equation*}
This contradicts that $w\in\Zlm_{\rbtab}$ and hence
we have $w{\rbtab}\in \TabRC$. 

Next, let us prove
 $\Zlm_{\rbtab}\rbtab\supseteq\TabRC$.
We will show that $w\in\Zlm_{\rbtab}$ for all $w$ such that
$w{\rbtab}\in \TabRC$ by induction on $l(w)$.
If $l(w)=0$ then $R(w)=\emptyset$ and $w\in\Zlm_{\rbtab}$.
Let $k\in\Z_{\geq 1}$ and
suppose that  the statement is true
for all $w$ with $l(w)<k$.

Take $w\in\affW$ such that $w{\rbtab}\in\TabRC$ and $l(w)=k$.
Take $x\in\affW$ and $i\in[0,n-1]$ 
such that $w=s_ix$ and $l(w)=l(x)+1$.
By Lemma~\ref{lem;wBRC}, we have 
$x{\rbtab}\in\TabRC$.
By the induction hypothesis, we have $x\in\Zlm_{\rbtab}$. 
Since  $R(w)=R(x)\sqcup\{x^{-1}(\al\ch_i)\}$,
it is enough to prove 
$$\sigma:=\brac{\wt_{\rbtab}}{x^{-1}(\alch_i)}=
\cont(x^{-1}(i))-\cont(x^{-1}(i+1))
\neq\pm 1.
$$

We put $T=xT_0$ in the rest of the proof.

Suppose 
$\sigma=1$.
Put $(a,b)=T^{-1}(i)$.
Then $T^{-1}(i+1)=(a+j+1,b+j)$ for some $j\in\Z$.
If $j<0$ then 
we have
 $(a,b-1)\in\alsm$ and
 $i+1=T(a+j+1,b+j)\leq T(a,b-1)<T(a,b)=i$.
This is a contradiction.
If $j>0$ then $(a+1,b)\in\alsm$ and
$i+1>T(a+1,b)>i$.
This is a contradiction too.
Therefore we must have $j=0$ and 
hence $T^{-1}(i+1)=(a+1,b)$.
But then 
we have 
$$w{\rbtab}(a,b)=s_iT(a,b)=i+1>i=s_iT(a+1,b)=w{\rbtab}(a+1,b)$$ and 
this contradicts  the assumption $wT_0\in\TabRC$.
Therefore $\sigma\neq 1$.

Suppose 
$\sigma=-1$. Put $(a,b)=T^{-1}(i)$.
Then similar argument as above
implies that $T^{-1}(i+1)=(a,b+1)$.
This yields a contradiction too. 

Therefore we have $w\in\Zlm_{\rbtab}$ and proved 
$\TabRC=\Zlm_{\rbtab}\rbtab$.

\medskip\noindent
({\it Step 2})
By Step1,
for each $T\in\TabRC$,
there exists 
$w_T\in\Zlm_{\rbtab}$ such that $T=w_T\rbtab$.

First we will show that $zw_T^{-1}\in\Zlm_T$ for all $z\in\Zlm_{\rbtab}$.

Assume that  $zw_T^{-1}\notin\Zlm_T$ for some $z\in\Zlm_{\rbtab}$. 
Then there exists $\al\in R(zw_T^{-1})$ such that $\brac{\wt_T}{\alch}=\pm 1$.

If  $\al\in w_T R^+$, then putting $\beta=w_T^{-1}(\al)$,
we have $\beta\in R(z)$ and 
$\brac{\wt_{\rbtab}}{\beta\ch}=
{\brac{\wt_{w_T \rbtab}}{w_T(\beta\ch)}   }
=\pm 1$.
This contradicts the choice $z\in\Zlm_{\rbtab}$. 

If $\al\notin w_T R^+$, then putting $\beta=-w_T^{-1}(\al)$,
we have $\beta\in R(w_T)$ and $\brac{\wt_{\rbtab}}{\beta\ch}=\pm 1$.
This contradicts the choice $w_T\in\Zlm_{\rbtab}$. 

Therefore $zw_T^{-1}\in \Zlm_T$ for all $z\in\Zlm_{\rbtab}$.
Similarly, one can show that $zw_T\in\Zlm_{\rbtab}$ for
all $z\in\Zlm_T$.
Hence the correspondence $z\mapsto zw_T^{-1}$ gives a bijection
from $\Zlm_{\rbtab}$ to $\Zlm_T$, whose
inverse is given by $z\mapsto zw_T$.
Therefore, we have
$$zT=z w_T \rbtab\in\TabRC\Leftrightarrow 
z w_T\in \Zlm_{\rbtab}\Leftrightarrow 
z \in\Zlm_T.
$$
\qed
\subsection{Content of standard tableaux}
\label{subsection-content.standard}

Let
$n\in\Z_{\geq2}$, $m\in\Z_{\geq 1}$ and $\l\in\Z_{\geq0}$.

Let  $\lmm\in\Icalml$ and $T\in\TabRC$.

Put $\kappa=\l+\m$ and  $F=\con_T$.
Then it is easy to check 
that $F$ satisfies the following:

  \smallskip
\noindent
${\rm{(C1)}}$ $F(i+n)=F(i)-\kappa$ for all $i\in \Z$.

\smallskip
\noindent
${\rm{(C2)}}$ 
For any $\p\in \Z$ and 
$i,j\in F^{-1}(\p)$ such that $i<j$ and 
$[i,j]\cap F^{-1}(\p)=\{i,j\}$,
there exist unique
$k_-\in F^{-1}(\p-1)$ and unique $k_+\in F^{-1}(\p+1)$ such that
$i<k_-<j$ and $i<k_+<j$ respectively.

\smallskip 
Notice that condition ${\rm{ (C1)} }$ implies that
 $\sharp F^{-1}(\p)$ is finite for all $\p\in\Z$.

Conversely, suppose that a map $F:\Z\to\Z$ 
satisfying the conditions ${\rm{(C1) 
(C2)}}$ 
is given.
Then it  can be seen that
$F$ is a content associated with
 some standard tableau
on some periodic skew diagram, 
as in the following proposition.
\begin{proposition}\label{pr;stcontent}
Let $n\in\Z_{\geq2}$ and $\kappa\in\Z_{\geq 1}$.
Suppose that the map $F:\Z\to\Z$ 
satisfies the conditions ${\rm{(C1) 
(C2)}}$ above.
Then there exist $m\in[1,\kappa]$, $\lmm\in\Icala_{\m,\kappa-\m}$
and $T\in\TabRC$ such that
$F=\con_T$.
\end{proposition}
\noindent{\it Proof.}
({\it Step 1})
For $\p\in F(\Z)\overset{def}{=}\{F(i)\in\Z\mid i\in\Z\}$,
put $d_p=\sharp F^{-1}(\p)$, which
 ${\rm{(C1)} }$ implies is finite.
 Let $i_\p^{(1)},i_\p^{(2)},\dots,i_\p^{(d_\p)}$ be
the integers 
such that 
 $i_\p^{(1)}<i_\p^{(2)}<\dots<i_\p^{(d_\p)}$ and
 $$F^{-1}(\p) =\{i_\p^{(1)},i_\p^{(2)},\dots,i_\p^{(d_\p)}\}.$$

It follows from the condition (C1) that
$
 d_\p=d_{\p-\kappa}
$
and
 \begin{equation}\label{eq;i+n}
i_{\p-\kappa}^\bri{1}=i_{\p}^\bri{1}+n,\
i_{\p-\kappa}^\bri{2}=i_{\p}^\bri{2}+n,\dots,
i_{\p-\kappa}^\bri{d_{\p-\kappa}}=
i_{\p}^\bri{d_\p}+n\quad    
 \end{equation}
for all $\p\in F(\Z)$.

The following statement follows easily from
the condition (C2) and an induction argument (on $j$):

\smallskip
\noindent
{\it Claim.
Let $p\in F(\Z)$.

\smallskip\noindent
$\rm{(i)}$ If $p+1\in F(\Z)$ and $i_p^{(1)}<i_{p+1}^{(1)}$,
then 
 $d_p-d_{p+1}=0$ or $1$, and it holds that
 \begin{align*}
 i_p^{(j)}&<i_{p+1}^{(j)}\ \ \ \,(j\in[1,d_{p+1}]),\\ 
 i_p^{(j)}&>i_{p+1}^{(j-1)}\ \ (j\in[2,d_p]).
 \end{align*}
$\rm{(ii)}$ If $p+1\in F(\Z)$
and $i_p^{(1)}>i_{p+1}^{(1)}$,
then
$d_p-d_{p+1}=0$ or $-1$, and it holds that
\begin{align*}
i_p^{(j)}&>i_{p+1}^{(j)}\ \ (j\in[1,d_{p}]),\\
i_p^{(j-1)}&<i_{p+1}^{(j)}\ \ (j\in[2,d_{p+1}]).  
\end{align*}
$\rm{(iii)}$ If $p+1\notin F(\Z)$
then $d_p=1$.}

\medskip\noindent
({\it Step 2})
Fix $\p_0\in F(\Z)$ and $r\in \Z$.
We will define a subset $\diagram=\diagram_{\p_0,r}$ of $\Z^2$ as follows:

For $\p\in F(\Z)$ define $\nxt$ as the minimum number
in $F(\Z)\cap \Z_{>\p}$.

There exists a unique sequence
$\{(a_\p^\bri1,b_\p^\bri1)\}_{\p\in F(\Z)}$ in $\Z^2$ 
satisfying the initial condition
\begin{equation}\label{eq;initial}
(a_{\p_0}^\bri1,b_{\p_0}^\bri1)=(r,\p_0+r)
\end{equation}
and  the recursion relation 
\begin{equation}\label{eq;condition_for_ab1}
  (a_{\nxt}^\bri1,b_{\nxt}^\bri1)=
\begin{cases}
(a_\p^\bri1,b_\p^\bri1+1)  &\hbox{if  } \nxt=\p+1
\hbox{ and }i_\p^\bri1<i_{\p+1}^\bri1,\hfill \hbox{(i)}\\
(a_\p^\bri1-1,b_\p^\bri1)  &\hbox{if  } \nxt=\p+1
\hbox{ and }i_\p^\bri1>i_{\p+1}^\bri1,\hfill \hbox{(ii)}\\
(a_\p^\bri1-1,b_\p^\bri1+\nxt-\p-1) \!\!&\hbox{if  }\nxt>\p+1.\hfill 
\hbox{(iii)}
\end{cases}
\end{equation}
\def\bx{\framebox(10,10)}
\begin{center}
\setlength{\unitlength}{.08cm}
\begin{picture}(60,45)(-20,-15)
\put(-65,15){(i) $i_\p^\bri1<i_{\p+1}^\bri1$}
\put(-60,-5){\bx{$i_\p^{(1)}$}}
\put(-50,-5){\bx{$i_{\p+1}^{(1)}$}}

\put(-10,15){(ii) $i_\p^\bri1>i_{\p+1}^\bri1$}
\put(0,-10){\bx{$i_\p^{(1)}$}}
 \put(0,0){\bx{$i_{\p+1}^{(1)}$}}
\put(50,15){(iii) $\nxt>\p+1$}
\put(50,-10){\bx{$i_\p^{(1)}$}}
\put(70,0){\bx{$i_\nxt^{(1)}$}}
\end{picture}
\end{center}
Put 
\begin{equation}
  \label{eq;condition_for_abj}
  \ab{\p}{j}=(a_\p^\bri1+j-1,b_\p^\bri1+j-1)
\quad (\p\in F(\Z),\ j\in[2,d_\p]),
\end{equation}
and put
\begin{equation}
  \diagram=\left\{\left.\ab{\p}{j}\in\Z^2
\right|
 \p\in F(\Z),\, j\in[1,d_\p]
\right\}.
\end{equation}
Note that  $\ab{p}{1}$ will be the most
northwest box in $\alsm$ on the diagonal with content $p$.

\medskip\noindent
({\it Step 3})
Now, we will check
that the set $ \diagram$
satisfies the
condition (D1)(D2)(D3) in Definition~\ref{def;periodicdiagram}.

\smallskip
\noindent
{\it Check} (D1):
For $\p\in F(\Z)$, put 
$$\l_\p=\sharp\{s\in[\p,\p+\kappa-1]\cap F(\Z)\mid
\tilde{s}=s+1\text{ and }i_s^{(1)}<i_{s+1}^{(1)}\},$$
and put $\m_\p=\kappa-\l_\p$.
Then $\l_p,\m_p\in[0,\kappa]$ and
$m_\p=a_{\p}^\bri1-a_{\p+\kappa}^\bri1$ by
\eqref{eq;condition_for_ab1}.
Moreover it follows 
from \eqref{eq;i+n} 
that the number
$\l_p$ is independent of $\p$, and so is $\m=\m_\p$.

Since $b_\p^{(1)}-a_\p^{(1)}=\p$,
we have
$\ab{\p-\kappa}{1}=(a_{\p}^\bri{1}+\m,
b_\p^\bri{1}-\kappa+\m),$
and hence
\begin{equation}
  \label{eq;abkappa}
\ab{\p-\kappa}{j}=(a_{\p}^\bri{j}+\m,
b_\p^\bri{j}-\kappa+\m)
 \end{equation}
for all $j\in[1,d_\p]$.
Therefore $\diagram$ satisfies the condition (D1)
with $\peri:=(\m,-\kappa+\m)$.

\noindent
{\it Check} (D2):
Put $E=
\{\ab{\p}{j}\mid \p\in[1,\kappa]\cap F(\Z),\ j\in[1,d_\p]\}$.
Then $E$
gives
a fundamental domain  of the action of $\Z \peri$
on  $\diagram$,
and the 
set $E$ is in one to one correspondence with the set 
$F^{-1}([1,\kappa])$, by the definition of $d_\p$.
Hence we have 
$\sharp E=\sharp F^{-1}([1,\kappa])=n$ by (C1) 
and thus the condition (D2)
is checked.

\noindent
{\it Check} (D3):
Note that 
Claim above implies that
\begin{equation}
  \label{eq;(11)}
  (a,b),(a+1,b+1)\in\diagram\Rightarrow
(a+1,b),(a,b+1)\in\diagram.
\end{equation}
Suppose that the condition (D3) does not hold.
Then there exist
$(a,b)\in\diagram$ and $(i,j)\in\Z_{\geq 0}\times\Z_{\geq 0}\setminus
\{(0,0),(1,0),(0,1)\}$
for which it holds that
\begin{align*}
&(a+i',b+j')\in\diagram\ \Leftrightarrow\
(i',j')=(0,0)\text{ or }(i,j).
\end{align*}
Fix such $(a,b)$ and $(i,j)$.

First, suppose that $j-i=0$.
Then by \eqref{eq;condition_for_abj},
 $i(=j)$ must be $1$.
This implies that
$(a+1,b),(a,b+1)\in\diagram$.
This is  a contradiction and
hence $i-j\neq 0$.

Next, suppose that $j-i>0$.
Let $(a,b)=\ab{s}{r}$ and $(a+i,b+j)=\ab{\p}{k}$.
Note that $\p-s=j-i> 0$.

If $k=1$ then 
we have $a_\p^{(1)}-a_s^{(r)}=i\geq 0$. 
On the other hand, it follows from the definition 
\eqref{eq;condition_for_ab1}\eqref{eq;condition_for_abj}
of $\{\ab{p}{j}\}_{\p\in F(\Z),j\in[1,d_\p]}$
that
$a_s^{(r)}\geq a_s^{(1)}\geq a_{\p}^{(1)}$
and the equalities hold only if
$r=1$ and 
$s,s+1,\dots,\p-1 \in F(\Z)$
and $i_{s}^{(1)}<i_{s+1}^{(1)}<\dots<i_{\p}^{(1)}$.
This
implies that $i=0$ and 
$\ab{s+j'}{1}=(a,b+j')\in\diagram$ for all $j'\in[0,j]$.
This is a contradiction.

If $k\neq 1$ then $\ab{\p}{k-1}=(a+i-1,b+j-1)\in\diagram$
and hence $(a+i,b+j-1)\in\diagram$ by \eqref{eq;(11)}.
This is a contradiction since
$(i,j-1)\neq (0,0),(i,j)$.

By similar argument, a contradiction is derived when $j-i<0$.

This means that (D3) holds for $\diagram$, and hence
$\diagram=\dset^{n}_{(\m,-\kappa+\m)}$.

Let us show that $\diagram$ contains no empty rows.
It is clear 
that
$\diagram$ contains empty rows only if $\m=0$, that
implies 
$F(\Z)=\Z$ and 
$i_p^{(1)}<i_{p+1}^{(1)}$ for all $p\in\Z$ by
\eqref{eq;condition_for_ab1}.
But then it follows that $i_{p-\kappa}^{(1)}<i_p^{(1)}$
and this contradicts \eqref{eq;i+n}.
Therefore we have $\diagram\in\dset^{*n}_{(\m,-\kappa+\m)}$,
or equivalently, 
 $\diagram=\alsm$ for some  $\lmm\in\Icalaml$.

\medskip
\noindent
({\it Step 4})
Define the map
 $T:\diagram\to\Z$ by  $T\ab{\p}{j}=i_{\p}^{(j)}$. 
Obviously, we have
 $F=C\circ T^{-1}$.
It follows from 
\eqref{eq;i+n} and \eqref{eq;abkappa} that
 $T$ is a tableau on $\diagram$.
Moreover, Claim in Step 1 implies
that $T$ is row-column increasing, namely,  $T\in\TabRC$. 
This completes the proof.
\qed

\medskip
For $\m\in\Z_{\geq 1}$, define an 
automorphism
$\omega_\m$ of $\Z^m$  by
\begin{equation}
  \omega_\m\cdot\lm=(\lm_\m+\l+1,\lm_1+1,\lm_2+1,\dots,\lm_{\m-1}+1),
\end{equation}
for $\lm=(\lm_1,\lm_2,\dots,\lm_\m)\in\Z^m$.
Let ${\bra \omega_\m\ket}$ denote
the free group generated by $\omega_\m$,
and
let  ${\bra \omega_\m\ket}$ act on $\Z^m\times \Z^m$ by 
$\omega_\m\cdot\lmm
=(\omega_\m\cdot\lm,\omega_\m\cdot\mu)$
for $\lmm\in\Z^m\times\Z^m$.
Note that
 ${\bra \omega_\m\ket}$
 preserves the subsets $\Icalml$ and $\Icalaml$ of $\Z^\m\times\Z^\m$.
\begin{proposition}\label{pr;diagonal_shift}
Let $\m,\m'\in[1,n]$ and $\l,\l'\in\Z_{\geq0}$.
  Let $\lmm\in\Icalaml$ and $(\eta,\nu)\in\Icala_{\m',\l'}$.
The following are equivalent$:$

\smallskip\noindent
$\rm{(a)}$ $\con^{\alsm}_T=\con^{\widehat{\eta/\nu}}_S$ for some $T\in\TabRC$ and 
$S\in\tab^{\rm{RC}}(\widehat{\eta/\nu})$,

\smallskip\noindent
$\rm{(b)}$
$\m=\m'$, $\l=\l'$ and $\alsm=\widehat{\eta/\nu}+(r,r)$ for
some $r\in\Z$.

\smallskip\noindent
$\rm{(c)}$ 
$\m=\m'$, $\l=\l'$ and $(\eta,\nu)=\omega_\m^r\cdot(\lm,\mu)$
 for some $r\in\Z$.
\end{proposition}
\noindent
{\it Proof.}
First we will prove (a)$\Leftrightarrow$(b).

It is easy to see that (b) implies (a).
To see that (a) implies (b), recall 
the proof of Proposition~\ref{pr;stcontent},
where
the relations 
\eqref{eq;condition_for_ab1}\eqref{eq;condition_for_abj}
together with the initial condition \eqref{eq;initial}
determine the periodic skew diagram  $\diagram_{\p_0,r}
=\{
\ab{\p}{j}\mid \p\in F(\Z), j\in[1,d_\p]\}$ 
and its period
uniquely for each $\p_0\in F(\Z)$ and $r\in\Z$.

Note that
\begin{equation}
  \label{eq;gamma}
\diagram_{\p_0,r'}=\diagram_{\p_0,r}+(r'-r,r'-r)
\end{equation}
for $r,r'\in\Z$.

Put $F=\con_T$.
As in the proof of Proposition~\ref{pr;stcontent},
we put
$d_p=\sharp F^{-1}(\p)$ for $\p\in F(\Z)$, and
let $i_\p^{(1)}<i_\p^{(2)}<\dots<i_\p^{(d_\p)}$ be
the integers 
such that 
 $F^{-1}(\p) =\{i_\p^{(1)},i_\p^{(2)},\dots,i_\p^{(d_\p)}\}.$

Put $\ab{\p}{j}=T^{-1}(i_\p^{(j)})$.
Then it is easy to see that the sequence 
$\{\ab{\p}{j}\}_{\p\in F(\Z), j\in[1,d_\p]}$ 
satisfies the relations 
\eqref{eq;condition_for_ab1}\eqref{eq;condition_for_abj}. 
Therefore we have 
$$\alsm=
\{T^{-1}(i_\p^{(j)})\mid \p\in F(\Z), j\in[1,d_\p]\}
=\diagram_{\p_0,r}$$ for some $\p_0\in F(\Z)$ and 
$r\in\Z$.
Similarly, we have 
$\widehat{\eta/\nu}=\diagram_{\p_0,r'}$ for some 
$r'\in \Z$ (with the same $\p_0\in F(\Z)$).
Now, it follows from \eqref{eq;gamma} that (a) implies (b).

The equivalence (b)$\Leftrightarrow$(c) follows from 
Proposition~\ref{pr;affine_diagram} and 
the formula
\omitted{
\begin{equation}
\widehat{\omega_\m^r\cdot\lm\, /\, \omega_\m^r\cdot\mu}=
\alsm-(r,r)\quad (r\in\Z),
\end{equation}
}
\begin{equation}
\begin{picture}(130,0)
\put(0,0){\mbox{$\omega_\m^r\cdot\lm\, /\, \omega_\m^r\cdot\mu$}}
\put(80,0){\mbox{$= \alsm-(r,r)\quad (r\in\Z)$,}}
\put(0,10){\line(5,1){35}}
\put(70,10){\line(-5,1){35}}
\end{picture}
\end{equation}
that is verified by a simple calculation.
\qed
\section{Representations of the double affine Hecke algebra}
\label{section-rep.DAHA}
Let $\k$ denote a field whose characteristic is not equal to $2$.
\subsection{Double affine Hecke algebra of type $A$}
\label{subsection-def.DAHA}
%
Let $q\in\k$.

The double affine Hecke algebra was introduced by 
Cherednik~\cite{Ch;unification,Ch;double}.
\begin{definition}
Let $n\in\Z_{\geq 2}$.  

\smallskip\noindent
$\rm{(i)}$
The {\it double affine Hecke algebra} $\daff H_n(q)$
of $GL_n$ is 
the unital
 associative algebra over $\k$ defined by the following
generators and relations:
\begin{alignat*}{2}
&{generators}\
:\ & &
 t_0,t_1,\dots, t_{n-1},\pi^{\pm 1},x_1^{\pm 1},x_2^{\pm 1},
\dots,x_n^{\pm 1},\cent^{\pm 1}.\\
&relations\ for\ n\geq 3\ 
&:\ &
(t_i-q)(t_i+1)=0~(i\in[0,n-1]),\\
 & & & t_it_j t_i=t_j t_i t_j~(j\equiv i\pm 1 \hbox{ mod }n),\\
 & & & t_it_j=t_j t_i~(j\not\equiv i\pm 1 \hbox{ mod }n),\\
 & & & \pi\pi^{-1}=\pi^{-1}\pi=1,\\
 & & &\pi t_i\pi^{-1}=t_{i+1}~(i\in[0,n-2]),\ \pi t_{n-1}\pi^{-1}=t_0,\\
 & & &x_i x_i^{-1}=x_i^{-1}x_i=1~~(i\in[1,n]),\\
 & & &x_ix_j=x_jx_i~(i,j\in[1,n]),\\
 & & &t_ix_it_i=qx_{i+1}~(i\in[1,n-1]),\ 
t_0x_{n}t_0=\cent^{-1}qx_1\\
 & & &t_ix_j=x_jt_i~(j\not\equiv i,i+1\hbox{ mod }n),\\
 & & &\pi x_i\pi^{-1}=x_{i+1}~(i\in[1,n-1]),\ 
\pi x_n\pi^{-1}=\cent^{-1}x_1,\\
 & & &\cent\cent^{-1}=\cent^{-1}\cent=1,\quad
 \cent^{\pm 1}h=h\cent^{\pm 1}~(h\in \ddot H_n(q)).\\
&relations\ for\ n=2\ 
&:\ &
 (t_i-q)(t_i+1)=0~(i\in[0,1]),\\
 & & &\pi\pi^{-1}=\pi^{-1}\pi=1,\quad
 \pi t_0\pi^{-1}=t_{1},\ \pi t_{1}\pi^{-1}=t_0,\\
 & & &x_i x_i^{-1}=x_i^{-1}x_i=1~~(i\in[1,2]),\quad
 x_1x_2=x_2x_1,\\
 & & &t_1x_1t_1=qx_{2},\quad t_0x_{2}t_0=\cent^{-1}qx_1\\
 & & &\pi x_1\pi^{-1}=x_{2},\  \pi x_2\pi^{-1}=\cent^{-1}x_1,\\
 & & &\cent\cent^{-1}=\cent^{-1}\cent=1,\quad
 \cent^{\pm 1}h=h\cent^{\pm 1}~(h\in \ddot H_2(q)).
\end{alignat*}
\noindent
$\rm{(ii)}$ 
Define the {\it affine Hecke algebra} $\dot H_n(q)$
of $GL_n$ as the subalgebra of $\daff H_n(q)$
generated by $\{t_0,t_1,\dots,t_{n-1},\pi^{\pm 1}\}$.
\end{definition}
\begin{remark}
  It is known that the subalgebra  of $\daff H_n(q)$
generated by
$$\{t_1,t_2,\dots,t_{n-1},x_1^{\pm 1},x_2^{\pm 1},\dots,x_n^{\pm 1}\}$$ is
also isomorphic to $\dot H_n(q)$.
\end{remark}
For $\nu=\sum_{i=1}^n\nu_i\e_i+ \nu_c c^*
\in \aff P$,
put
$$x^{\nu}=x_1^{\nu_1}x_2^{\nu_2}\dots x_n^{\nu_n}\xi^{\nu_c}.$$

Let $\X$ denote the commutative group
$\{x^\nu \mid \nu\in\aff P\}\subseteq \daff H_n(q)$.
The group algebra $\k[\X]
=\k[x_1^{\pm1},x_2^{\pm1},\dots,x_n^{\pm1},\xi^{\pm1}]$
is a commutative subalgebra 
of $\daff H_n(q)$.

For $w\in\aff W$ 
with a reduced expression
 $w=\pi^r s_{i_1}s_{i_2}\cdots s_{i_k}$, put
$$t_w=\pi^r t_{i_1}t_{i_2}\cdots t_{i_k}.$$
Then $t_w$ does not depend on the choice of the reduced expression,
and $\{t_w\}_{w\in\aff W}$ forms a basis of the affine
Hecke algebra $\aff H_n(q)\subset \daff H_n(q)$.

It is easy to see that 
$\{t_w x^\nu\}_{w\in \aff W,\nu\in \aff P}
$ and $\{x^\nu t_w\}_{w\in \aff W,\nu\in \aff P}$
respectively form basis of $\daff H_n(q)$.
In particular, we have 
\begin{proposition}\label{pr;PBW}
$\daff H_n(q)=\aff H_n(q)\k[\X]=\k[\X]\aff H_n(q)$.
\end{proposition}
Define an element $\phi_i$ of $\daff H_n(q)$ by
\begin{align}
  \phi_i&=t_i
\left(
1-x^{\al_i}
\right)
+1-q\quad (i\in[0,n-1]).
\end{align}
By direct calculations, we have the following:
\begin{lemma}\label{lem;Phi_relation}
The following holds in $\daff H_n(q):$
  \begin{align*}
&\phi_i\phi_j=\phi_j\phi_i
\quad (i,j\in[0,n-1],\ 
j\not\equiv i\pm 1\hbox{ mod }n),\\
&\phi_i\phi_j\phi_i=\phi_j\phi_i\phi_j
\quad (i,j\in[0,n-1],\ 
j\equiv i\pm 1\hbox{ mod }n),\\
&\phi_i^2
={(1-qx^{\al_i})(1-qx^{-\al_i})}
\quad(i\in[0,n-1]).  
  \end{align*}
\end{lemma}
For $w\in\aff W$ with a 
reduced expression $w=\pi^r s_{i_1}s_{i_2}\cdots s_{i_k}$, 
put
$$\phi_w={\pi^r}\phi_{i_1}\phi_{i_2}\cdots\phi_{i_k}.$$
Then $\phi_w$ does not depend on the choice of
the reduced expression by Lemma~\ref{lem;Phi_relation}.
For an $\daff H_n(q)$-module $M$,
the element  $\phi_w\in\daff H_n(q)$ is 
regarded as a linear  operator on $M$,
and 
$\phi_w$ is called an {\it intertwining operator}.

The following formula follows easily:
\begin{lemma}
$\phi_w x^\nu=x^{w(\nu)}\phi_w$
for any $w\in\aff W$ and  $\nu\in\aff P$. 
\end{lemma}
%
%
%
\begin{lemma}  \label{lem;inter_in_t}
For $w\in\aff W$,
\begin{equation}
\phi_w =t_w\prod_{\al\in R(w)}(1-x^\al)+
\sum_{y\in\aff W,y\prec w}t_y f_{y}    
\end{equation}
for some $f_{y}\in\k[\X]$. 
\end{lemma}
\proof
Follows from the expression \eqref{eq;R(w)} of $R(w)$
and induction on $l(w)$.
\qed

\medskip

Let $\X^*$ denote
 the set of characters of $\X$:
$$\X^*=\Hom_{group}(\X,\mathrm{GL}_1(\k)).$$
Consider the correspondence
$\aff P\to\X^*$
 which maps 
$\zeta
\in \aff P$ to the character
$q^\zeta\in\X^*$ defined by
$$q^\zeta(x_i)= q^\brac{\zeta}{\ech_i}\ 
(i\in[1,n]),\quad
 q^\zeta(\xi)=q^\brac{\zeta}{c},
$$
or equivalently, defined by
$q^\zeta(x^\nu)= q^\brac{\zeta}{\nu\ch}\  (\nu\in\aff P).$
Through this correspondence,
$\aff P$ is identified with the subset
$$\{\chi\in\X^*\mid \chi(x^\nu)\in q^\Z\ (\forall \nu\in\aff P)
\}$$
 of
$\X^*$, where $q^\Z=\{q^r\mid r\in\Z\}$.

For an $\daff H_n(q)$-module
 $M$ and $\zeta
\in \aff P$, 
define the weight space $M_{\zeta}$ 
and the generalized weight space $M_{\zeta}^{\gen}$
of weight $\zeta$ 
with respect to the action of $\k[\X]$ by 
\begin{align*}
M_{\zeta}&=\left\{ v\in M \left|\
(x^\nu-q^{\brac{\zeta}{\nu\ch}}) v=0\
\hbox{ for any }  
\nu\in \aff P\right.\right\},\\
M_{\zeta}^{\gen}&=
\bigcup_{k\geq 1}
\left\{ v\in M \left| \
(x^\nu-q^{\brac{\zeta}{\nu\ch}})^k v=0\hbox{ for any }  \nu\in\aff P
\right.\right\}.
\end{align*}
For an $\daff H_n(q)$-module $M$, an element $\zeta\in \aff P$ is called
a weight of $M$ if $M_\zeta\neq 0$,
and an element $v\in M_\zeta$ (resp.\;$M_\zeta^\gen$)
 is called a weight vector (resp.\;generalized weight vector)
of weight $\zeta$.

The following statement can be verified by direct calculations.
\begin{proposition}\label{pr;invertible}
  Let $M$ be an $\daff H_n(q)$-module.
Let $\zeta\in \aff P$ and  $v\in M_\zeta$.
Then the following holds for all $w\in\aff W:$

\smallskip
\noindent
{\rm{(i)}}
$\phi_w M_\zeta\subseteq M_{w(\zeta)}$ and
$\phi_w M^{\gen}_\zeta\subseteq M^{\gen}_{w(\zeta)}.$

\smallskip
\noindent
{\rm{(ii)}} 
$\phi_{w^{-1}}\phi_w v=
\prod_{\al\in R(w)}
{(1-q^{1+\bra \zeta\mid \al\ch\ket })
(1-q^{1-\bra \zeta\mid \al\ch\ket})}v.
$
\end{proposition}
For $\zeta\in \aff P$,
put 
\begin{equation}
  \label{eq;Zzeta}
\aff\ZZ_{\zeta}=\{w\in \aff W
\mid \brac{\zeta}{\alch}\notin\{-1,1\}\hbox{ for all }\al\in R(w)\}.
\end{equation}
Note that $\aff\ZZ_{\wt_T}=\Zlm_T$ for
$\lmm\in\Icalml$ and $T\in\TabRC$.

As a direct consequence of 
Proposition~\ref{pr;invertible}, we have the following:
\begin{proposition}\label{pr;ZZ}
Suppose that $q$ is not a root of $1$.
Let $M$ be an $\daff H_n(q)$-module and  $\zeta\in \aff P$.
For $w\in\aff\ZZ_\zeta$, 
the map
$$\phi_w:M_\zeta\to M_{w(\zeta)}$$
is a linear isomorphism.
\end{proposition}
\subsection{$\X$-semisimple modules}
 \label{section-Xssl}
\begin{remark}
 Througout Section\;\ref{section-Xssl}, the lemmas and propositions are still
true and require almost no modification of their statements or proofs, even
if $\kappa$ is not an integer or if $q$ is a root of unity.
However, we impose these restrictions so that the combinatorics
developed in Section\;\ref{section-periodic.skew}
describe the structure of the $\X$-semisimple modules.
When we relax the condition $\kappa \in \Z$ but still require
$q$ generic, one can extend the combinatorial description with
appropriate reformulation.
\end{remark}

Fix $n\in\Z_{\geq 2}$.
Let $q\in\k$ and suppose that $q$ is not a root of $1$.

Fix $\kappa\in\Z$ and
put $P_\kappa=P+ \kappa c^*=
\{\zeta\in\aff P\mid \brac{\zeta}{c}=\kappa\}$.
\begin{definition}
Define $\Cat$ as the set
consisting of those $\daff H_n(q)$-modules $M$ 
which is finitely generated and admits a decomposition
$$M=\bigoplus_{\zeta\in P_\kappa} M_\zeta$$
 with $\dim M_\zeta<\infty$ for all $\zeta\in P_\kappa$.
\end{definition}

We say that a module $M \in \Cat$ is {\it $\X$-semisimple}.
%
%

In the following, we will see 
 some general properties of 
$\daff H_n(q)$-modules in $\Cat$.
The results and argument used in the proofs
are essentially the same as those
for the affine Hecke algebra
(See e.g. \cite{Ram1}). 
\begin{lemma}\label{lem;weightnonzero}
   Let $M\in\Cat$.
Let $i\in[0,n-1]$ and let $\zeta\in P_\kappa$ be
such that
$\brac{\zeta}{\alch_i}= 0$.
Then $M_\zeta=\{0\}$.
\end{lemma}
\noindent
{\it Proof.}
Suppose that there exists $v\in M_\zeta\setminus\{0\}$.
Then we have
\begin{align*}
& (x^{\al_i}-1)t_i v=2(1-q)v\neq 0,\\
&  (x^{\al_i}-1)^2 t_i v=0.
\end{align*}
This implies
$t_iv\in M_\zeta^\gen\setminus M_\zeta$, which
contradicts the assumption $M=\oplus_{\zeta \in P_\kappa} M_\zeta$.
\qed
\begin{lemma}\label{lem;span} 
   Let $L$ be an irreducible $\daff H_n(q)$-module which belongs to
$\Cat$.
Let $v$ be a non-zero weight vector of $L$.
Then $L=\sum_{w\in \aff W}\k\phi_w v$.
\end{lemma}
\proof
Put $N=\sum_{w\in \aff W}\k\phi_w v\subseteq L$.
Since $L=\sum_{w\in \aff W}\k t_w v$
by Proposition~\ref{pr;PBW}, it is enough
to prove that 
$t_wv\in N$ for all $w\in \aff W$.
We proceed by induction on $l(w)$.

It is clear
 that $t_wv\in N$  for $w$ of length zero.

Let $k\in\Z_{\geq 1}$ and 
suppose that $t_wv\in N$ for all $w\in \aff W$ with $l(w)<k$.
Take $w\in \aff W$ with a reduced expression 
$w=\pi^r s_{i_1}s_{i_2}\dots s_{i_k}$ (and hence $l(w)=k$).
By Lemma~\ref{lem;inter_in_t},
we have $\phi_w v=
\sum_{x\in\aff W,\;x\preceq w}g_{wx} t_xv$ with some 
coefficients
$g_{wx}\in\k$.

If $g_{ww}\neq 0$ then 
$t_wv=g_{ww}^{-1}(\phi_w v -\sum_{x\prec w}
g_{wx}t_xv)\in N$.

Suppose $g_{ww}= 0$.
By Lemma~\ref{lem;inter_in_t}, this means
\begin{equation}
  \label{eq;f_ww}
\prod_{\al\in R(w)}(1-q^{\brac{\zeta}{\alch}})=0,
\end{equation}
where $\zeta\in P_\kappa$ is the weight of $v$.
Hence, there exists $p\in[1,k]$ such that
\begin{align}\label{eq;fyy}
 &\prod_{\al\in R(y)}(1-q^{\brac{\zeta}{\alch}})\neq 0,\quad
\prod_{\al\in R(s_{i_p}y)}(1-q^{\brac{\zeta}{\alch}})= 0,
\end{align}
where $y=s_{i_{p+1}}s_{i_{p+2}}\dots s_{i_k}$.
This implies $\brac{\zeta}{y^{-1}(\alch_{i_p})}=
\brac{y(\zeta)}{\alch_{i_p}}=0$.
By Lemma~\ref{lem;weightnonzero},
 we have $L_{y(\zeta)}=0$ and hence 
$\phi_y v=0.$

Let 
 $\phi_y v=\sum_{x\in\aff W,\;x\preceq y}g_{yx} t_xv$
with $g_{yx}\in \k$.
Multiplying $\pi^r t_{i_1}t_{i_2}\dots t_{i_p}$
to the equality $\phi_y v=0$, we have
\begin{equation}
  \label{eq;t_wv}
g_{yy}t_wv=-\sum_{x\prec y}g_{yx}
\pi^r t_{i_1}t_{i_2}\dots t_{i_p}
t_xv.
  \end{equation}
Note that $g_{yy}=
\prod_{\al\in R(y)}(1-q^{\brac{\zeta}{\alch}})
\neq 0$ by \eqref{eq;fyy}, and
it is easy to verify that
the right hand side of \eqref{eq;t_wv}
is in $N$ using the induction hypothesis.
Therefore $t_wv\in N$.
\qed

\smallskip
For $\zeta\in P_\kappa$, let 
 $\aff W[\zeta]$ denote the stabilizer of $\zeta$:
$$\aff W[\zeta]=\{w\in\aff W\mid w(\zeta)=\zeta\}.$$
\begin{lemma}\label{lem;phizero}
Let $L$ be an irreducible $\daff H_n(q)$-module which belongs to
$\Cat$.
Let $\zeta$ be a weight of $L$ and let $v\in L_\zeta$.
Then $\phi_wv=0$ for all $w\in\aff W[\zeta]\setminus \{1\}$. 
\end{lemma}
\proof
Let $w\in\aff W[\zeta]\setminus \{1\}$ with 
a reduced expression
 $w=\pi^r s_{i_1}s_{i_2}\dots s_{i_k}$.

Put
$\aff R[\zeta]=\{\al\in\aff R\mid \brac{\zeta}{\alch}=0\}.$

Then, $\aff R[\zeta]$ is a subroot system of $\aff R$ and
$\aff W[\zeta]$ is the corresponding Coxeter group.
Moreover it follows that
a system of positive (resp.\;negative) roots 
is given by
$\aff R[\zeta]\cap \aff R^+$ (resp.\;$\aff R[\zeta]\cap \aff R^-$).

Therefore for $w\in\aff W[\zeta]\setminus \{1\}$, there exists
a reflection $s_\al\ (\al\in \aff R[\zeta]\cap \aff R^+)$
such that $w(\al)\in \aff R[\zeta]\cap \aff R^-\subseteq \aff R^-$.

Now, Lemma~\ref{lem;coxeter}-(ii)
implies
that there exists $p\in[0,n-1]$ such that 
$ws_\al=\pi^r s_{i_1}s_{i_2}\dots s_{i_{p-1}}s_{i_{p+1}}\dots s_{i_k}$.
Putting $y=s_{i_{p+1}}s_{i_{p+2}}\dots s_{i_k}$, we have
$\brac{\zeta}{y^{-1}(\alch_{i_p})}=\brac{y(\zeta)}{\alch_{i_p}}=0$.
Lemma~\ref{lem;weightnonzero} implies that $L_{y(\zeta)}=0$
and hence
 $\phi_{w}v=\phi_{\pi^r s_{i_1}s_{i_2}\dots s_{i_{p}}}
\phi_y v=0$.
\qed
\begin{proposition}\label{pr;onedim}
  Let $L$ be an irreducible $\daff H_n(q)$-module which belongs to
$\Cat$.
Then $\dim L_\zeta\leq 1$ 
for all $\zeta\in P_\kappa$.
\end{proposition}
\proof
Directly follows from Lemma~\ref{lem;span}
and Lemma~\ref{lem;phizero}.
\qed
\begin{lemma}\label{lem;pm1}
  Let $L$ be an irreducible $\daff H_n(q)$-module which belongs to
$\Cat$.
Let $\zeta$ be a weight of $L$ and let $i\in[0,n-1]$
such that $\brac{\zeta}{\alch_i}\in\{-1,1\}$. Then $\phi_iv=0$
for $v\in L_\zeta$.
\end{lemma}
\proof
Suppose $\brac{\zeta}{\alch_i}=\pm 1$ and let $v\in L_\zeta\setminus\{0\}$.

Suppose $\phi_i v\neq 0$.
Put $\aff W'=\{w\in\aff W\mid ws_i\in \aff W[\zeta]\}$.
Then it follows from Lemma~\ref{lem;span} 
that
\begin{equation}\label{eq;sumW'}
  \sum_{w\in \aff W'}a_w\phi_w\phi_i v=v
\end{equation}
for some $\{a_w\in\k\}_{w\in\aff W'}$.

For $w\in\aff W'$ such that
 $l(w s_i)<l(w)$, we have
$\phi_w=\phi_{w s_i}\phi_i$.
Proposition~\ref{pr;invertible}-(ii) implies that
$$\phi_w\phi_i v=\phi_{s_iw}\phi_i^2 v
=\phi_{s_iw}
(1-q^{1+\brac{\zeta}{\alch_i}})(1-q^{1-\brac{\zeta}{\alch_i}})
v$$
and it is $0$ as $\brac{\zeta}{\alch_i}=\pm 1$.

For $w\in\aff W'$ such that $l(ws_i)>l(w)$, 
we have $\phi_w\phi_i v=\phi_{ws_i}v=0$
by Lemma~\ref{lem;phizero}.

Therefore the left hand side of \eqref{eq;sumW'} is $0$
and this is a contradiction.
\qed
\subsection{Representations associated with periodic skew diagrams}
\label{ss;construction}
%
In the rest of this paper, we always assume that $q$
is not a root of $1$.

Let $n\in\Z_{\geq 2}, \ \m\in\Z_{\geq 1}$ and $\l\in\Z_{\geq 0}$.

For $\lmm\in\Icalml$, set 
\begin{equation}
  \V\lmm=\bigoplus_{T\in\TabRC}\k v_T.
\end{equation}
Define linear operators  $\tilde{x}_i~(i\in[1,n])$, 
$\tilde{\pi}$ and
 $\tilde{t}_i~(i\in[0,n-1])$ 
on $\V\lmm$ by
\begin{align}
\tilde{x}_i v_T&=q^{\con_T(i)}v_T,\\
\tilde{\pi} v_T&=v_{\pi T},\\
\tilde{t}_i v_T&=
\begin{cases}
\frac{1-q^{1+\tau_i}}{1-q^{\tau_i}}v_{s_i T}
-\frac{1-q}{1-q^{\tau_i}}v_T\quad
&\text{ if }s_iT\in\TabRC,\\
-\frac{1-q}{1-q^{\tau_i}}v_T\quad
&\text{ if }s_iT\notin\TabRC,
\end{cases}
\label{eq;Taction}
\end{align}
where 
$$
\tau_i
=\con_T(i)-\con_T(i+1)
=\brac{\wt_T}{\alch_i}
\quad (i\in[0,n-1]).
$$
The following lemma is easy
and ensures that the operator
 $\tilde t_i$ is well-defined:
\begin{lemma}
$\con_T(i)-\con_T(i+1)\neq 0$ for any $i\in[0,n-1]$ and $T\in\TabRC$.
\end{lemma}
\begin{theorem}
Let $\lmm\in\Icalml$.
  There exists an algebra homomorphism 
$\homlm:\daff H_n(q)\to\End_\k(\V\lmm)$
such that
\begin{alignat*}{2}
&\homlm(t_i)=\tilde t_i\ (i\in[0,n-1]),\ 
& &\homlm(\pi)= \tilde \pi,\\
&\homlm(x_i)=\tilde x_i\ (i\in[1,n]),\ 
& &\homlm(\cent)= q^{\l+m}.
\end{alignat*}
\end{theorem}
\noindent{\it Proof.}
The defining relations  of
$\daff H_n(q)$ can be verified by direct calculations
(See \cite{Ram1}, for a sample of calculation for the affine Hecke algebra).
\qed

Note that  the $\daff H_n(q)$-module
$\V\lmm$ for $\lmm\in\Icalml$  belongs to $\Cat$ with $\kappa=\l+\m$.

%
\begin{theorem}\label{th;irredandmfree}
Let $\lmm\in\Icalml$.

\smallskip\noindent
$\rm{(i)}$
$\V\lmm=\bigoplus_{T\in\TabRC}\V\lmm_{\wt_T}$, and
$\V\lmm_{\wt_T}=\k v_T$ for all $T\in\TabRC$.

\smallskip
\noindent
$\rm{(ii)}$ The $\daff H_n(q)$-module $\V\lmm$ is irreducible.
\end{theorem}
\noindent{\it Proof.}
(i)
Follows directly from Proposition~\ref{pr;CT=CS}.

(ii) Let $N$ be a non-zero submodule of $\V\lmm$.
Since $N$ contains at least one weight vector,
we can assume that $v_T\in N$ for some $T\in\TabRC$.

Let $S\in\TabRC$.
By Theorem~\ref{th;Zlm=TabRC2}, 
there exists $w_S\in\Zlm_T$ 
such that $S=w_S T$. 
Put $\tilde v_S=\phi_{w_S}v_T\in N$.
Since the intertwining operator
$$\phi_{w_S}:\V\lmm_{\wt_T}\to \V\lmm_{w_S(\wt_T)}
=\V\lmm_{\wt_S}$$
 is a linear isomorphism
by Proposition~\ref{pr;ZZ},
we have $\tilde v_S\in  \V\lmm_{\wt_S}\setminus\{0\}$.

Now, it follows from (i) that 
$\oplus_{S\in\TabRC}\k \tilde v_S=\V\lmm$.
Therefore we have
$N\supseteq \V\lmm$ and hence $N=\V\lmm$.
Therefore $\V\lmm$ is irreducible.
\qed
\subsection{Classification of $\X$-semisimple modules}
\label{subsection-classification}
Fix $n\in\Z_{\geq 2}$ and $\kappa\in\Z_{\geq1}$.
Let $q\in\k$ and suppose that $q$ is not a root of $1$.

Our next and final purpose is to show that the modules
$\V\lmm$ we constructed in Section\;\ref{ss;construction}
exhausts all irreducible modules in $\Cat$.
\begin{lemma}\label{lem;klemma}
Let $L$ be an irreducible $\daff H_n(q)$-module 
which belongs to $\Cat$.
For any weight $\zeta\in P_\kappa$ of $L$
and $i,j\in\Z$ such that $i<j$ and
$$\brac{\zeta}{\alch_{ij}}=0,$$
there exist $k_+\in [i+1,j-1]$ and  
$k_-\in[i+1,j-1]$ such that
$$\brac{\zeta}{\alch_{i\,k_+}}=-1\text{ 
and } 
\brac{\zeta}{\alch_{i\,k_-}}=1$$ respectively.
\end{lemma}
%
\noindent{\it Proof.}
We proceed by induction on $j-i$.

For any weight $\zeta$ of $L$ and $i\in\Z$,
we have $\brac{\zeta}{\alch_i}\neq 0$ by Lemma~\ref{lem;weightnonzero}.
Therefore we have nothing to prove when $j-i=1$.

Let $r>1$ and assume that the statement holds when $j-i<r$.

In order to complete the induction step,
it is enough to prove the existence of $k_\pm$ for
a weight $\zeta$ of $L$ and $i,j\in \Z$ such that $j-i=r$ 
and
\begin{equation}
  \label{eq;j-imin}
  \{k\in[i+1,j-1]\mid \brac{\zeta}{\alch_{i\,k}}=0\}=\emptyset.
\end{equation}
Fix a non-zero weight vector $v\in L_\zeta$.

\noindent
Case (i): Suppose
$\brac{\zeta}{\alch_i}=\pm 1$ and $\brac{\zeta}{\alch_{j-1}}= \pm 1$.

Then the statement holds with $k_\pm=j-1$ and $k_\mp=i+1.$

\noindent
Case (ii): Suppose
$\brac{\zeta}{\alch_i}= -1$ 
and $\brac{\zeta}{\alch_{j-1}}= 1$.

Then we have $\brac{\zeta}{\alch_{i+1\,j-1}}=0$.
If $i+1\neq j-1$ then there exist
$k'_-\in[i+1,j-1]$ such that $
\brac{\zeta}{\alch_{i+1\,k'_-}}= 1$,
and hence $\brac{\zeta}{\alch_{i\,k'_-}}
=\brac{\zeta}{\alch_{i}}+
\brac{\zeta}{\alch_{i+1\, k'_-}}=0$.
This contradicts the choice \eqref{eq;j-imin} of $i,j$.
Therefore we have $j-i=2$.
This case, 
we have $\brac{\zeta}{\alch_i}=-1$ and 
$\brac{\zeta}{\alch_{i+1}}=1$.
Hence Lemma~\ref{lem;pm1} implies that
$\phi_i v=0$ and $\phi_{i+1}v=0$, 
which gives
 $t_iv=qv$ and
$t_{i+1}v=-v$
respectively.
But then we have 
$$-q^2v=t_it_{i+1}t_iv=t_{i+1}t_it_{i+1}v=qv,$$
and this is a contradiction as $q$ is not a root of $1$.
Therefore 
this case is not possible.

\noindent
Case (iii):  Suppose  
$\brac{\zeta}{\alch_i}= 1$ and $\brac{\zeta}{\alch_{j-1}}=-1$.

A similar argument as in (ii) implies that this case
is not possible.

\noindent
Case (iv):
Suppose $\brac{\zeta}{\alch_i}\neq \pm1$.

Then $\phi_iv\neq 0$ by Proposition~\ref{pr;invertible} and
hence $s_i(\zeta)$ is a weight of $L$. 
By $\brac{s_i(\zeta)}{\alch_{i+1\,j}}=0$,
the induction hypothesis implies that
there exists 
{$k_\pm\in[i+2,j-1]$ } 
such that $
\brac{\zeta}{\alch_{i\,k_\pm}}=
\brac{s_i(\zeta)}{\alch_{i+1\,k_\pm}}=\mp 1$.
Hence the statement holds.

\noindent
Case (v): Suppose $\brac{\zeta}{\alch_{j-1}}\neq \pm1$.

Then $\phi_{j-1}v\neq 0$ and a similar argument 
as in (iv)
implies that there exists 
{$k_\pm\in[i+1,j-2]$ }
 such that $\brac{\zeta}{\alch_{i\,k_\pm}}=\mp1$.

This completes the proof 
\qed
%
\begin{theorem}\label{th;completerep}
Let $n\in\Z_{\geq 2}$
and $\kappa\in\Z_{\geq 1}.$
  Let $L$ be an irreducible $\daff H_n(q)$-module
which belongs to $\Cat$.
Then there exist $\m\in[1,\kappa]$
and $\lmm\in\Icala_{\m,\kappa-\m}$ such that
$L\cong \V\lmm$.
\end{theorem}
\noindent
{\it Proof.}
\noindent
({\it Step 1})
Let $\zeta\in P_\kappa$ be a weight of $L$.

Define $F_\zeta:\Z\to \Z$ by $F_\zeta(i)=\brac{\zeta}{\ech_i}$
$(i\in\Z)$.

It is easy to see that $F_\zeta$ satisfies  condition
(C1) 
 in Proposition~\ref{pr;stcontent},
and  the existence of $k_\pm$ in
condition (C2) 
follows from
Lemma~\ref{lem;klemma}.
Note that the uniqueness of $k_\pm$ in (C2) 
follows automatically
from the condition that
$[i,j]\cap F_\zeta^{-1}(\p)=\{i,j\}$, setting $p = F_\zeta(i)$ here.
Suppose, without loss of generality, that there were another choice of
$k_\pm'$, with $k_\pm < k_\pm'$.  It follows 
 $\brac{\zeta}{\alch_{k_\pm \,k_\pm'}}=0$, and applying Lemma~\ref{lem;klemma}
here gives the existence of an $i'$ between $k_\pm$ and $k_\pm'$ and hence
$i < i' < j$ with 
 $\brac{\zeta}{\alch_{i\,i'}}=0$. This gives 
$i' \in F_\zeta^{-1}(\p)$, a contradiction.

Therefore, Proposition~\ref{pr;stcontent} implies
that there exist $\m\in[1,n]$, 
 $T\in\Tab$ and $\lmm\in\Icala_{\m,\kappa-\m}$
such that
$F_\zeta=\con_T$, or equivalently,
$\zeta=\wt_T$.


\noindent
({\it Step 2})
Recall that $\aff\ZZ_\zeta=\Zlm_T$.

Take $u\in L_\zeta\setminus\{0\}$.
For each $w\in\Zlm_T$,
put 
\begin{align}
 \sigma_w&=\prod_{\al\in R(w)}(1-q^{1+\brac{\zeta}{\alch}})
\end{align}
and
\begin{align}
u_w&=\sigma_w^{-1} \phi_{w} u. 
\end{align}
Here, note that $\sigma_w\neq 0$ and $u_w\neq 0$
for all $w\in\Zlm_T$
by Proposition~\ref{pr;ZZ}.

Put
$N=\sum_{w\in\aff\ZZ_\zeta}\k \phi_{w}u
=\sum_{w\in\Zlm_T}\k u_w\subseteq L.$
Since
$u_w\in L_{\wt_{wT}}$
and each weight space is linearly independent
by Proposition~\ref{pr;CT=CS},
we have
$N=\bigoplus_{w\in\Zlm_T}\k u_w.$

By Theorem~\ref{th;Zlm=TabRC2}, 
one can define $w_S\in\Zlm_T$ by $S=w_ST$ for all $S\in\TabRC$, 
and
define 
a linear map $\bijec:\V\lmm\to L$
by $\bijec(v_S)=u_{w_S}$  $(S\in\TabRC)$.
It is obvious that 
$\bijec$ is injective and its image is $N$.

Let us see that $\bijec$ is an $\daff H_n(q)$-homomorphism.

Let $w\in\Zlm_T$.
Let $i\in[0,n-1]$ be such that $l(s_iw)<l(w)$.
Then we have $s_iw\in\Zlm_T$ and
$\sigma_w=(1-q^{1+\brac{\zeta}{(s_iw)^{-1}(\alch_i)}})
\sigma_{s_iw}$.
Therefore
\begin{equation}
\begin{split}\label{eq;Phiiuw1}
\phi_i u_w
&=
\sigma_w^{-1}\phi_i  \phi_w u
=\sigma_w^{-1} \phi_i^2 \phi_{s_iw}u\\
&=
( 1-q^{ 1+\brac{s_i w(\zeta) }{\alch_i} })
( 1-q^{ 1-\brac{s_i w(\zeta) }{\alch_i} })
\sigma_{w}^{-1}\phi_{s_iw}u\\
&=
( 1-q^{ 1+\brac{w(\zeta)}{\alch_i} })
u_{ s_iw }.
\end{split}
\end{equation}
Let $i\in [0,n-1]$ be such that $l(s_iw)>l(w)$.
If  $s_iw\not\in\Zlm_T$  then $\brac{\zeta}{w^{-1}(\alch_i)}=\pm 1$
and hence $\phi_i u_w=0$ by Lemma~\ref{lem;pm1}.
If  $s_iw\in\Zlm_T$ then 
 we have
 $\sigma_{s_iw}=(1-q^{1-\brac{\zeta}{w^{-1}(\alch_i)}})\sigma_w$
and 
\begin{align*}
\phi_i u_w&=
\sigma_w^{-1} \phi_i\phi_w u
=\sigma_w^{-1} \phi_{s_iw}u\\
&=(1-q^{1+\brac{w(\zeta)}{\alch_i} })u_{ s_iw }.
\end{align*}
Therefore, in both cases, we have
\begin{equation*}\label{eq;Phiiuw}
 \phi_i u_w=  \begin{cases}
(1-q^{1+\brac{w(\zeta)}{\alch_i} })u_{ s_iw }
\quad & (s_iw\in\Zlm_T),\\
0     & (s_iw\not\in\Zlm_T).
\end{cases}
\end{equation*}
This implies
$$\bijec({t_i}v_S)=
t_i\bijec(v_S)\ \ (i\in [0,n-1],\; S\in\TabRC).$$
Moreover, 
it is easy to see that 
\begin{align*}
\bijec({x}_i v_S)&=x_i\bijec(v_S)\ (i\in[1,n]),\quad
\bijec({\xi} v_S)=\xi\bijec(v_S),\quad
\bijec({\pi} v_S)=\pi\bijec(v_S)
\end{align*}
for all $S\in\TabRC$.
Therefore $\bijec$ is 
an $\daff H_n(q)$-homomorphism and it gives an 
isomorphism $\V\lmm\cong N$ of $\daff H_n(q)$-modules.
Since $L$ is irreducible, we have $L=N\cong\V\lmm$.
\qed
%
%
\begin{corollary}
  Let $L$ be an irreducible $\daff H_n(q)$-module
which belongs to $\Cat$.
Let $v\in L$ be a non-zero
weight vector of weight $\zeta\in \aff P_\kappa$.
Then 
$$L=\bigoplus_{w\in \aff\ZZ_{\zeta}} \k \phi_w v,$$  
and $\phi_w v\neq 0$ for all $w\in\aff\ZZ_{\zeta}$.
\end{corollary}
\begin{theorem}\label{th;distinct}
Let
 $m,m'\in\Z_{\geq 1}$ and $\l,\l'\in\Z_{\geq 0}$.
Let $\lmm\in\Icala_{\m,\l}$ and
 $(\eta,\nu)\in\Icala_{\m',\l'}$.
Then the following are equivalent$:$

\smallskip\noindent
$\rm{(a)}$ $\V\lmm\cong \V(\eta,\nu)$. 

\smallskip\noindent
$\rm{(b)}$ 
$\m=\m'$, $\l=\l'$ and
$\alsm=\widehat{\eta/\nu}+(r,r)$ for some $r\in\Z$.

\smallskip\noindent
$\rm{(c)}$ 
$\m=\m'$, $\l=\l'$  and $(\eta,\nu)=\omega_\m^r\cdot\lmm$
for some $r\in\Z$.
\end{theorem}
\noindent
{\it Proof.}
Follows from (Step1) in the proof of Theorem~\ref{th;completerep}
and Proposition~\ref{pr;diagonal_shift}.
\qed

\medskip
Let $\Irr\Cat$ denote the set of isomorphism classes
of all simple modules in $\Cat$.
Combining Theorem~\ref{th;completerep} and Theorem~\ref{th;distinct},
we obtain the following classification theorem,
which is announced in \cite{Ch;fourier} in more general situation.
\begin{corollary}[cf.\cite{Ch;fourier}]\label{cor;clasification}
Let $n\in\Z_{\geq 2}$ and $\kappa\in\Z_{\geq 1}$.
The correspondences $\lmm\mapsto \alsm$ and 
 $\lmm\mapsto \V\lmm$ induce
the following bijections respectively$:$
$$\bigsqcup_{\m\in[1,\kappa]}\dset^{*n}_{(\m,-\kappa+\m)}/\Z(1,1)
\overset{\sim}{\leftarrow}
\bigsqcup_{\m\in[1,\kappa]}\left(\Icala_{\m,\kappa-\m}/\bra\omega_\m\ket\right)
\overset{\sim}{\rightarrow} \Irr\Cat.$$
\end{corollary}
\begin{remark}
We gave a direct and combinatorial proof for
Theorem~\ref{th;completerep}, Theorem~\ref{th;distinct}
and Corollary~\ref{cor;clasification}
based on the tableaux theory on periodic skew diagrams.

An alternative approach to prove
these results
is to use the result in \cite{Va,Su}, where
the classification of irreducible modules over $\daff H_n(q)$
of a more general class is obtained.
Actually, it is easy to see
that the $\daff H_n(q)$-module $\V\lmm$
coincides with the unique simple quotient $\daff L\lmm$
 of the induced module $\daff M\lmm$ with the notation in \cite{Su}.
\end{remark}
\begin{remark}
It is easy to derive the corresponding results
for the degenerate affine Hecke algebra by 
a parallel argument.
\end{remark}
\begin{remark}
  There exists an algebra involution
$\iota:\daff H_n(q)\to\daff H_n(q)$ such that
\begin{align*}
  &\iota(t_i)=qt_i^{-1}~(i\in[0,n-1]),\
\iota(\pi)=\pi,\\ 
&\iota(x_i)=x_i^{-1}~(i\in[1,n]),
\ \iota(\cent)=\cent^{-1}.
\end{align*}
The composition $\homlm\circ\iota:\daff H_n(q)\to \V\lmm$
gives an $\daff H_n(q)$-module structure on $\V\lmm$
on which $\xi$ acts as a scalar $q^{-\l-\m}$.
We let $\V^\iota\lmm$ denote this   $\daff H_n(q)$-module.
The correspondence $\lmm\mapsto \V^\iota\lmm$ 
 induces a bijection 
$$\bigsqcup_{\m\in[1,\kappa]}
\left(\Icala_{\m,\kappa-\m}/\bra\omega_\m\ket\right)
\to\Irr\cat_{-\kappa}(\daff H_n(q))$$
for all $\kappa\in\Z_{\geq 1}$.
\end{remark}

\end{document}